\documentclass[a4paper]{article}

\usepackage[plainpages=false, colorlinks=true, 
            linkcolor=black, urlcolor=black, citecolor=black]{hyperref}
\usepackage{geometry}
\usepackage{tabularx}
\usepackage{tikz}
\usepackage{overpic}
\usepackage{amsmath, amssymb, amsthm,  amsfonts, mathrsfs, cite}
\usepackage{multirow}
\usepackage{color}
\usepackage{booktabs}
\usepackage{tikz}
\usepackage{bm}
\usepackage{hyperref}
\usepackage{scrextend}
\usepackage{arydshln}
\usepackage {url}
\usepackage[font=scriptsize, labelfont=bf]{caption}

\usepackage{graphicx} 
\usepackage{epstopdf} 
\usepackage{algpseudocode}
\usepackage{graphicx}
\usepackage{array,  tabularx}
\usepackage{algorithm, algcompatible}
\usepackage{booktabs} 
\usepackage{paralist} 
\usepackage{verbatim} 
\usepackage{subfig} 

\usepackage{xcolor, marginnote, enumitem}
\usepackage[numbers]{natbib}
\numberwithin{equation}{section}
\theoremstyle{plain}
\newtheorem{lemma}{Lemma}
\newtheorem{theorem}{Theorem}
\newtheorem{proposition}{Proposition}

\theoremstyle{remark}
\newtheorem{remark}{Remark}

\theoremstyle{definition}
\newtheorem{definition}{Definition}

\DeclareMathOperator{\gph}{gph}
\DeclareMathOperator{\dist}{dist}

\DeclareMathOperator{\proj}{proj}
\DeclareMathOperator{\prox}{prox}
\DeclareMathOperator{\dom}{dom}

\DeclareMathOperator*{\argmin}{arg\, min}

\newcommand{\bb}{\mathbb}

\newcommand{\lookUp}[1]{}

\newcommand{\bitem}{\begin{itemize}}
	\newcommand{\eitem}{\end{itemize}}

\newcommand{\bpm}{\begin{pmatrix}}

\DeclareMathAlphabet{\mathbfit}{OML}{cmm}{b}{it}

\usepackage{xspace}
\usepackage{bold-extra}
\usepackage[most]{tcolorbox}

\colorlet{texcscolor}{blue!50!black}
\colorlet{texemcolor}{red!70!black}
\colorlet{texpreamble}{red!70!black}
\colorlet{codebackground}{black!25!white!25}

\date{}

\begin{document}

\title{A preconditioned difference of convex functions algorithm with extrapolation and line search}

\author{Ran Zhang \thanks{School of Mathematics, 
		Renmin University of China,  China  \  \href{mailto:shenxinhua@ruc.edu.cn}{zhangran227@ruc.edu.cn}}
\and   Hongpeng Sun\thanks{School of Mathematics,
Renmin University of China, China \  \href{mailto:hpsun@amss.ac.cn}{hpsun@amss.ac.cn} }
}

\maketitle

\begin{abstract}
This paper proposes a novel proximal difference-of-convex (DC) algorithm enhanced with extrapolation and aggressive non-monotone line search for solving non-convex optimization problems. We introduce an adaptive conservative update strategy of the extrapolation parameter determined by a computationally efficient non-monotone line search. The core of our algorithm is to unite the update of the extrapolation parameter with the step size of the non-monotone line search interactively. The global convergence of the two proposed algorithms is established through the Kurdyka-\L ojasiewicz properties, ensuring convergence within a preconditioned framework for linear equations. Numerical experiments on two general non-convex problems: SCAD-penalized binary classification and graph-based Ginzburg-Landau image segmentation models, demonstrate the proposed method’s high efficiency compared to existing DC algorithms in both convergence rate and solution accuracy.
\end{abstract}

\paragraph{MSCodes.}{65K10,
65F08,
49K35,
90C25,
90C26.
}


\section{Introduction}\label{sec:intro}

We aim to solve the following non-convex optimization problem
\begin{equation}\label{min_problem}
\min_{x\in X} \ E(x): = f( x ) +g( x) \,=\,f\left( x \right) +g_1\left( x \right) -g_2\left( x \right).
\end{equation}
Here $g_1\left( x \right)$ and $g_2\left( x \right)$ are proper closed convex functions, 
{$f$ is convex and} 
 $\nabla f\left( x \right) $  is Lipschitz continuous with constant $L$, and $X$ is a finite-dimensional Hilbert space. Such non-convex optimization formulations find broad applications in digital image processing, statistical machine learning, and related fields \cite{Apx,BF,Pangcui}. The proximal DC algorithm (DCA) with extrapolation (shortened as pDCAe henceforth) has emerged as a particularly efficient approach, successfully accelerating conventional proximal gradient methods   \cite{LeThi2018,LPH,Lethi2024}. This method combines the computational efficiency of proximal gradient algorithms with DC decomposition techniques, as detailed in \cite{Apx,Pangcui} for various non-convex optimizations. Our work builds upon these foundations while introducing critical enhancements to the algorithmic framework.

As shown in the survey \cite[Section 5]{Lethi2024}, the convergence speed of DCA can be accelerated through two distinct approaches: extrapolation techniques with traditional parameterization strategies, including FISTA-based extrapolation \cite[Algorithm 1]{ijcai2018p190} or FISTA-based restart mechanisms, e.g., \cite[Algorithm 5]{Lethi2024}, \cite[Algorithm (3.4)]{Wen2018} and Armijo line search methods based on decreasing energy \cite[Algorithm 2]{Artacho2018}, \cite{shen2024preconditioned}.

Compared to the existing works, our contributions are as follows. 
The main contributions are a concise {non-monotone line search mechanism \cite{ferreira2024boosted} and a new adaptive update strategy for the extrapolation parameter. The line search completely determines the extrapolation parameter.
Inspired by \cite{Artacho2018, ArtachoSIAM,ferreira2024boosted}, we employ the aggressive non-monotone line search, where the starting point for line search is generated by the preconditioned and extrapolated DCA, instead of using the previous iteration \cite{Artacho2018,ArtachoSIAM}.
The aggressive non-monotone line search \cite[Example 4]{ferreira2024boosted} can not only relax stringent descent requirements, but can also obtain an extra step. The proposed framework enhances pDCAe by integrating these two accelerating methods. Specifically, the non-monotone line search strategy provides global guidance for local extrapolation parameter updates. We refer to \cite[NLSA]{ZhangHongchao} \cite[Algorithm 1]{lu2019nonmonotone} \cite[Section 4]{Dai2002} for the traditional non-monotone line search strategy. Moreover, given the limited existing work that combines these approaches, our comparative analysis demonstrates distinct advantages over existing methods. Compared to adaptive linear search methods \cite[Algorithm 1]{yang2024proximal}, which require solving auxiliary subproblems, our approach eliminates computational overhead while preserving the core efficiency of pDCAe.

Furthermore, we introduce a line search framework for the case where $g_1$ is non-smooth, providing a global convergence guarantee that is not covered by the line search framework \cite{ArtachoSIAM}. In \cite{ArtachoSIAM}, only the case $g_2$ is non-smooth, and $g_1$ with Lipschitz gradient is discussed to ensure the existence of descent directions.

Numerical validations confirm the improvements. In the numerical experiments in Section \ref{sec:num}, we compare our algorithm with the algorithms that use only a single type of acceleration method. The experiments show that our algorithm has better performance. First, for iterative sequences near critical points, the integrated update strategy stabilizes convergence by controlling the extrapolation parameter $\beta$, effectively addressing the performance degradation observed in the original implementations of pDCAe. Second, the interactive combination of extrapolation and line search accelerates convergence compared to standard extrapolation and restart approaches. These empirical results, detailed in Section \ref{sec:num}, verify that our framework maintains the computational simplicity of pDCAe while achieving better convergence across various non-convex optimization benchmarks.

The paper is organized as follows. Section \ref{sec:intro} establishes notations and preliminary materials essential for the subsequent analysis. Section \ref{sec:extra-line} introduces two distinct DC decomposition algorithms incorporating a non-monotone line search strategy and a novel update strategy of the extrapolation parameter. Section \ref{sec:conver} rigorously analyzes the convergence properties of the proposed algorithms and examines a generalized computational framework. Finally, Section \ref{sec:num} validates the theoretical developments through numerical experiments, and Section \ref{sec:conclu} gives a conclusion.

\section{Notation and preliminaries}\label{sec:intro}
Throughout this paper, we use $X$ to denote the $n$-dimensional Hilbert space with inner product $\langle  \cdot, \cdot \rangle $ and Hilbert norm $\left\| \cdot \right\| $. Moreover, for a symmetric matrix $M\in \mathbb{R} ^{n\times n}$, we use $\mu _{\max}\left( M \right) $ and $\mu _{\min}\left( M \right) $ to denote its largest and smallest eigenvalues respectively.

We will touch on some necessary tools from the convex and variational analysis    \cite{Roc1}.
The {graph} of a multi-valued mapping $F : \mathbb{R}^n {\rightrightarrows} \mathbb{R}^m$ is defined by
\[
\gph F := \{(x,y) \in \mathbb{R}^n \times \mathbb{R}^m : y \in F(x)\}, 
\]
whose {domain} is defined by $ \dom F: = \{ x \ | \ F(x) \neq \emptyset \} $. 
Similarly the graph of an extended real-valued function $f :\mathbb{R}^n \to \mathbb{R} \cup \{+\infty\}$ is defined by 
\[
\gph f := \{(x,s) \in \mathbb{R}^n\times \mathbb{R} : s = f(x)\}. 
\]
Let $h: \mathbb{R}^n \rightarrow \mathbb{R}\cup \{+\infty\}$ be a proper lower semicontinuous function. Denote $\dom h: = \{ x \in\mathbb{R}^n: \ h(x) < +\infty  \}$. 
For each $x\in\dom f$, 
the limiting-subdifferential of $h$ at $x\in \mathbb{R}^n$, written $\partial f$, is defined as follows    \cite{Roc1}, 
\begin{align}
\partial h(x): = \bigg\{ & \xi \in \mathbb{R}^n: \exists x_n \rightarrow x,  h(x_n) \rightarrow h(x),  \xi_n  \rightarrow \xi,  \lim_{y\rightarrow x } \inf_{y \neq x_n}\frac{h(y)-h(x_n) - \langle \xi_n, y-x_n \rangle }{|y-x_n|}\geq 0 \bigg\}.  \notag
\end{align}

We next recall the Kurdyka-\L ojasiewicz (KL) property, which is satisfied by a wide variety of functions, such as proper closed semi-algebraic functions. It plays an important role in the convergence analysis of many first-order methods.

Let $h: \mathbb{R}^n \rightarrow \mathbb{R}\cup \{+\infty\}$ be a proper lower semicontinuous function. Denote $\dom h: = \{ x \in\mathbb{R}^n: \ h(x) < +\infty  \}$. 
For each $x\in\dom f$, the limiting-subdifferential of $h$ at $x\in \mathbb{R}^n$, written $\partial f$ are well-defined as in    \cite[Definition 8.3(a)]{Roc1}. Additionally, if $h$ is continuously differentiable,
the subdifferential reduces to the usual gradient $\nabla h$.
We also need the following  Kurdyka-\L ojasiewicz (KL) property and KL exponent for the global and local convergence analysis. While the KL properties can help obtain the global convergence of iterative sequences,  the KL exponent can help provide a local convergence rate. 
\begin{definition}[KL property, KL function and KL exponent    \cite{ABS}] \label{def:KL}
	A proper closed function $h$ is said to satisfy the KL property at $\bar x \in \dom \partial h$ if there exists $\nu \in (0,+\infty]$, a neighborhood $\mathcal{O}$ of $\bar x$, and a continuous concave function $\psi: [0,\nu) \rightarrow [0,+\infty)$ with $\psi(0)=0$ such that
	\begin{itemize}
		\item [{{(i)}}] $\psi$ is continuous differentiable on $(0,\nu)$ with $\psi'>0$ over $(0,\nu)$;
		\item [{{(ii)}}] for any $x \in \mathcal{O}$ with $h(\bar x) < h(x) <h(\bar x) + \nu$, one has
		\begin{equation}\label{eq:kl:def}
		\psi'(h(x)-h(\bar x)) \dist(0,\partial h(x)) \geq 1. 
		\end{equation}
	\end{itemize}
	A proper closed function $h$ satisfying the KL property at all points in $\dom \partial h$ is called a KL function. Furthermore, for a proper closed function $h$ satisfying the KL property at $\bar x \in \dom \partial h$,  if $\psi$ in \eqref{eq:kl:def} can be chosen as $\psi(s) = cs^{1-\theta}$ for some $\theta \in [0,1)$ and $c>0$, i.e., there exist  $\bar c, \epsilon >0$ such that
	\begin{equation}\label{eq:KL:exponent:theta:exam}
	\dist(0,\partial h(x)) \geq \bar c (h(x)-h(\bar x))^{\theta}
	\end{equation}
 whenever $\|x -\bar x\| \leq \epsilon$ and $h(\bar x) < h(x) <h(\bar x) + \nu$, then we say that $h$ has the KL property at $\bar x$ with exponent $\theta$.
	If $h$ has the KL property with exponent $\theta$ at any $ \bar x \in \dom \partial h$, we call $h$ a KL function with exponent $\theta$. 
\end{definition}

We also need the level-boundedness. A function $F: \mathbb{R}^n \rightarrow [-\infty, +\infty]$ is level-bounded (see Definition 1.8    \cite{Roc1}) if lev$_{\leq \alpha}F: =\{ u: F(u) \leq \alpha\}$ is bounded (or possibly empty). Level boundedness is usually introduced to ensure the existence of a minimizer. Throughout this paper, we assume that $E(x)$ in \eqref{min_problem} is level-bounded. We also assume $E(x)$ in \eqref{min_problem} is bounded below, i.e., $\inf_{x\in X}E(x) > -\infty$.

\section{Preconditioned DCA with extrapolation determined  by line search}\label{sec:extra-line}
Now let us turn to preconditioned DCA with extrapolation for \eqref{min_problem}. By solving the subproblem involving the difference of convex functions by incorporating a proximal term to update $\left\{ \bar{x}^n \right\}_n $, the optimization task for determining $\bar{x}^n$ is as follows
\begin{align}\label{subproblem}
\bar{x}^n=\mathop {\argmin} \limits_{y\in \mathbb{R} ^n}\left\{ \left< -\xi ^n,y \right> +\frac{1}{2}\left\| y-y^n \right\| _{M}^{2} + f\left( y \right) +g_1\left( y \right) \right\},
\end{align}
where $\xi^n \in \partial g_2\left( x^n \right)$, $\forall n \in N$, and $y^n=x^n+\beta_n\left( x^n-x^{n-1} \right)$.

Henceforth, we assume that the symmetric positive definite matrix $M$ serves as a preconditioning weight operator, where  $\langle a,b \rangle _M=\langle a, Mb \rangle$ for $\forall$ $a, b\in X$ and the induced norm as $\left\| x \right\| _{M}^{2}:=\left< Mx,x \right> $. 
Applying first-order optimality conditions to the minimization problem \eqref{subproblem} yields the critical equation for determining $\bar{x}^n$
\begin{equation}\label{eq:topre}
    \xi ^n \in  M\left( y-y^n \right) +\nabla f\left( y \right) +\partial g_1\left( y \right) .  
\end{equation}
We assume the above equation admits a unique solution, yielding the following update 
\begin{equation}\label{update_x_bar}
\bar{x}^n=\left( M+\nabla f+\partial g_1 \right) ^{-1}\left( My^n+\xi ^n \right). 
\end{equation}
{Here $\nabla $ denotes the gradient operator, while $\nabla f(y)\in \mathbb{R} ^n$ represents the evaluated gradient vector at point $y$.}

If $g_1+f$ is a quadratic function, the update $\left\{ \bar{x}^n \right\}_n $ in \eqref{update_x_bar} can be reformulated as the classical efficient preconditioned iteration, which can circumvent solving the linear equation even with mild accuracy, especially for a large-scale linear system. We will explain the preconditioning technique in Proposition \ref{pro:pre} later.


Within the widely employed gradient descent framework, \cite{o2015adaptive} pioneered two distinct restart strategies for the extrapolation parameter: a gradient-based scheme \cite{hinder2020generic,fercoq2019adaptive,Wen2018} and a function-value-based scheme. The former, taking \cite{Wen2018} for example, represents a computationally efficient implementation, while the latter, a function value-based restart scheme, i.e., restarting when $E\left( x^n \right) >E\left( x^{n-1} \right)$, which conceptually parallels the line search techniques proposed here, provides enhanced theoretical guarantees. We also refer to \cite{aujol2024parameter} for other restart techniques with the estimation of Lipschitz constants based on the FISTA framework. However, the function-value-based restart scheme usually needs to be restarted frequently due to the strict energy-decreasing requirement, which leads to an increased computational cost, as observed in our numerical experiments. 
Crucially, we address a critical failure mode: persistent line search failures that trigger excessive restarts, degrading the algorithm to basic DCA performance. To mitigate this, we directly introduce a novel update strategy that dynamically adjusts the extrapolation parameter based on the following line search. By defining $d^n = \bar x^n-x^n$, we proceed with an $\emph{aggressive}$ non-monotone Armijo-type line search involving parameters $\rho^k \lambda_{\max}>0$, $k = 1,2, \ldots$, $\rho \in (0,1)$, and $\lambda_{\max}, \omega, \eta>0$. The aim is to find the smallest integer $k$ satisfying \cite{ferreira2024boosted}  
\begin{equation}\label{eq:line:search}
E(\bar x^{n}+\rho^{k-1} \lambda_{\max} d^n) \leq {E}(\bar x^{n}) - \eta \rho^{k-1} \lambda_{\max} \|d^{n}\|_2^2 +\nu_{n}, \quad\nu _{n}:=\frac{\omega}{n+1}\left\| d^{n} \right\| ^2.
\end{equation}
Henceforth, denote $a(n):=k$ as the number of line searches for the $n$-th iteration and $N_{\max}$ as the maximal number of line searches. 
The final step size $\lambda_{n} := \rho^{k-1} \lambda_{\max}$ obtained from the line search in \eqref{eq:line:search} dictates the update of $\left\{ x^n \right\}_n $ as follows
\begin{equation}\label{eq:LSstep}
   x^{n+1} = \bar x^n + {\lambda}_{n} d^{n}.
\end{equation}
We propose the following method to update the extrapolation parameter $\left\{ \beta _{n+1} \right\}$ determined by line search \eqref{eq:line:search}. We assume that $b_1>0$ and $b_2 \ge 0$ are constants close to 0 henceforth.

\textbf{Line search determined extrapolation (shortened as LSDE henceforth):}
\begin{equation} \label{LSDE}
\begin{aligned} 
&\text{Case 1: }1\le a(n) \leq N_{\max}\text{, }\lambda_{n} := \rho^{k-1} \lambda_{\max}\text{, we set } \beta _{n+1}=\frac{1}{1+b_1+\lambda _{n}}, \\
&\text{Case 2: }a(n)>N_{\max}\text{, }\lambda_n=0 \text{, we set } \beta _{n+1}= b_2.
\end{aligned}    
\end{equation}
{The proposed LSDE that focuses on the global energy can provide global guidance for local extrapolation parameter updates.}

A critical distinction between our algorithm and conventional FISTA-based methods lies in parameterizing the extrapolation update $\left\{ \beta _{n+1} \right\}$. Established schemes-including FISTA with fixed or adaptive restart-initialize $\theta _{-1}=\theta _0=1$ and recursively compute $\beta_{n+1}$ as
\begin{equation}\label{convention_restart}
\beta _{n+1}=\frac{\theta _n-1}{\theta _{n+1}}, \quad \theta _{n+1}=\frac{1+\sqrt{1+4\theta _{n}^{2}}}{2},
\end{equation}
with parameter resets $\theta _{n}=\theta _{n+1}=1$ under specific conditions: fixed restart triggers resets periodically every $\bar{T}$ iterations, whereas adaptive restart activates when $\left< y^{n}-x^{n+1},x^{n+1}-x^{n} \right> >0$. 
\begin{proposition}\label{prop1}
According to \eqref{LSDE} for updating the extrapolation parameter $\left\{ \beta _{n+1} \right\}$, we assert that there exists a positive constant $C_{\lambda}$ such that $\frac{1}{\left( 1+\lambda _n \right) ^2}-\beta _{n+1}^{2}>C_{\lambda}$ for $\forall n$.
\end{proposition}
\begin{proof}
Recalling that $\frac{1}{\left( 1+\lambda _n \right) ^2}>\beta _{n+1}^{2}$, we then discuss this in two scenarios.

If $\lambda _n\ne 0$, known that $\beta _{n+1}=\frac{1}{1+b_1+\lambda _n}$, we then obtain
\begin{align*}
\frac{1}{\left( 1+\lambda _n \right) ^2}-\frac{1}{\left( 1+b_1+\lambda _n \right) ^2}=\frac{b_1\left( 1+b_1+2\lambda _n \right)}{\left( 1+\lambda _n \right) ^2\left( 1+b_1+\lambda _n \right) ^2}
\ge \frac{b_1\left( 1+b_1 \right)}{\left( 1+\lambda _{\max} \right) ^2\left( 1+b_1+\lambda _{\max} \right) ^2}\triangleq C_{\lambda}^{1}.
\end{align*}
If $\lambda _n=0$, we know that $\beta _{n+1}=0$, thus we have \begin{equation*}
\frac{1}{\left( 1+\lambda _n \right) ^2}-\beta _{n+1}^{2}\ge 1-\varepsilon _0>C_{\lambda}^{2}.    
\end{equation*}
Setting $C_{\lambda}=\min \left( C_{\lambda}^{1},C_{\lambda}^{2} \right) $,   we proved this conclusion.
\qed
\end{proof}

Specifically, we introduce two distinct DC splitting algorithms, both integrated with a non-monotone line search strategy and our proposed update strategy \eqref{LSDE} for solving the minimization problem \eqref{min_problem}; see Algorithm \ref{algorithm1} and Algorithm \ref{algorithm2} for more details. If putting $f(y)$ in the implicit part as $g_1$, we get Algorithm  \ref{algorithm1}. While putting $f(y)$ in the explicit form with linearization as $g_2$ \cite{Wen2018}, we obtain Algorithm  \ref{algorithm2}.



\begin{algorithm}  
\renewcommand{\algorithmicrequire}{\textbf{Input:}}
\renewcommand{\algorithmicensure}{\textbf{Output:}}
\caption{Algorithmic framework for 1-order extrapolation and new-proposed convex splitting method with non-monotone line search and preconditioning (shortened as \protect{\textbf{npDCAe$_{\text{nls}}$}})}\label{algorithm1} 
\begin{algorithmic}[1]
\State Choose $x^0\in \dom g_1,  \left\{ \beta _n \right\} \subset \left[ 0,1 \right) $ and $\mathop {\sup} \limits_{n}\beta _n<1$, $x^{-1}=x^0$, $\beta_0=0$.
\State Choose $\xi ^n\in \partial g_2\left( x^n \right), \forall n\in N$, $y^n=x^n+\beta _n\left( x^n-x^{n-1} \right)$. 
\State Solve the following subproblem 
\begin{equation}\label{eq:algorithmic:update1}
    \bar{x}^n=\mathop {\argmin} \limits_{y\in X }\left\{ \left< -\xi ^n,y \right> +\frac{1}{2}\left\| y-y^n \right\|^2 _M+f\left( y \right) +g_1\left( y \right) \right\}. 
\end{equation}
\State {Set $d^n=\bar{x}^n-x^n$. IF $d^n = 0$, STOP and RETURN $x^{n+1}$.  }
\State Set $\lambda _{\max}>0$, $\eta >0$, $\rho \in \left[ 0,1 \right) $, $N_{\max}\in \mathbb{N} $.
\State Set $\nu _n=\frac{\omega}{n+1}\left\| d^n \right\| ^2$,  $a\left( n \right) =1$, $\lambda_n=\lambda_{\max}$.
\State Start the non-monotone line search:
\begin{equation}\label{non-mono1}
\text{IF \quad} E\left( \bar{x}^n+\lambda _nd^n \right) >E\left( \bar{x}^n \right) -\eta \lambda _n\left\| d^n \right\| ^2+\nu _n   
\end{equation}
\State Do $\lambda _n= \rho \cdot\lambda_{n}$, set $a\left( n \right) =a\left( n \right) +1$; Otherwise, RETURN $x^{n+1} := \bar{x}^n + \lambda_nd^n$.
\State If $a\left( n \right) >N_{\max}$, STOP and RETURN  $x^{n+1} := \bar{x}^n$; Otherwise, go to Step 8.
\State Update $\beta_{n+1}$ as LSDE \eqref{LSDE}.
\end{algorithmic}
\end{algorithm}

\begin{algorithm}  
\renewcommand{\algorithmicrequire}{\textbf{Input:}}
\renewcommand{\algorithmicensure}{\textbf{Output:}}
\caption{Algorithmic framework for 1-order extrapolation and original convex splitting method with non-monotone line search and preconditioning (shortened as \protect{\textbf{pDCAe$_{\text{nls}}$}})}\label{algorithm2} 
\begin{algorithmic}[1]
\State Choose $x^0\in \dom g_1,  \left\{ \beta _n \right\} \subset \left[ 0,1 \right) $ and $\mathop {\sup} \limits_{n}\beta _n<1$, $x^{-1}=x^0$, $\beta_0=0$.
\State Choose $\xi ^n\in \partial g_2\left( x^n \right),\forall n\in N$, $y^n=x^n+\beta _n\left( x^n-x^{n-1} \right)$. 
\State Solve the following subproblem 
\begin{equation}\label{eq:algorithmic:update2}
\bar{x}^n=\mathop {\argmin} \limits_{y\in X}\left\{ \left< \nabla f\left( y^n \right) -\xi ^n,y \right> +\frac{L}{2}\left\| y-y^n \right\| _{M}^{2}+g_1\left( y \right) \right\}. 
\end{equation}
\State Set $d^n=\bar{x}^n-x^n$. IF $d^n = 0$, STOP and RETURN $x^{n+1}$.
\State Set $\lambda _{\max}>0$, $\eta >0$, $\rho\in \left[ 0,1 \right) $, $N_{\max}\in \mathbb{N} $.
\State Set $\nu _n=\frac{\omega}{n+1}\left\| d^n \right\| ^2$,  $a\left( n \right) =1$, $\lambda_n=\lambda_{\max}$.
\State Start the non-monotone line search:
\begin{equation}\label{non-mono2}
\text{IF \quad} E\left( \bar{x}^n+\lambda _nd^n \right) >E\left( \bar{x}^n \right) -\eta \lambda _n\left\| d^n \right\| ^2+\nu _n   
\end{equation}
\State Do $\lambda _n= \rho \cdot\lambda_{n}$, set $a\left( n \right) =a\left( n \right) +1$; Otherwise, RETURN $x^{n+1} := \bar{x}^n + \lambda_nd^n$.
\State If $a\left( n \right) >N_{\max}$, STOP and RETURN  $x^{n+1} := \bar{x}^n$; Otherwise, go to Step 8.
\State Update $\beta_{n+1}$ as LSDE \eqref{LSDE}.
\end{algorithmic}
\end{algorithm}




To better show the novelty of our algorithms, we collect some commonly used DC algorithms and their respective characteristics in the following Table \ref{table:Comparasion}.
\begin{table}[htbp!] 
 \caption{Commonly used DC algorithms and their respective frameworks.}
\label{table:Comparasion}
\begin{tabularx}{40em}
 {|*{6}{>{\centering\arraybackslash}X|}}
 \hline 
accelerated approach & extrapolation & line search & preconditioner & algorithm  \\ \hline
I & FISTA (restart) & No & No& pDCAe \cite{Wen2018} \\ \hline
II & No & monotone & No & BDCA \cite{Artacho2018} \\ \hline
III & No & non-monotone & No & nmBDCA \cite{ferreira2024boosted} \\ \hline
IV & FISTA (restart) & No & Yes & preDCAe \cite{DS}\\ \hline
V & LSDE & non-monotone & Yes & Our work\\ \hline
 \end{tabularx}
\end{table}
Before starting the proof of convergence, we shall give the strategy \cite{ferreira2024boosted} to choose the parameter sequence $\{\nu _n\}_n$. 
\begin{remark}\label{nu_n}
Let $\omega>0$ be a constant, and the sequence $\left\{\nu_n \right\}_n$ defined by $\nu _n=\frac{\omega}{n+1}\left\| d^n \right\| ^2$, for all $n\in N$.
Indeed, due to $\lim_{n\rightarrow \infty} \nu _n=\lim_{n\rightarrow \infty} \frac{\omega}{n+1}=0$, for every $\delta>0$, there exists $n_0 \in N$ such that $n\geq n_0$ implies $\frac{\omega}{n+1}\le \delta < \eta \cdot \lambda _{\max}\cdot \rho^{N_{\max}} $. We then have $\nu_n\le \delta\left\|d^n\right\|^2$. Similarly, we can define $\nu _n=\frac{\omega}{p\left( n \right)}\left\| d^n \right\| ^2$, where $p\left( n \right) $ is a monotonically increasing function concerning $n$, here the choice of $p\left( n \right) $ depends on the rate at which we want the non-monotone term to decline. 
\end{remark}
This brings great flexibility. At the beginning of iterations,   the line search \eqref{eq:line:search} with $\nu_n$ is indeed non-monotone. However, it eventually becomes monotone for sufficiently large $n$ since $\nu_n$ converges to $0$ while $N_{\max}$ is bounded.

Now, we end this section with the following proposition for addressing the preconditioning.
\begin{proposition}\label{pro:pre}
Suppose $(g_1+f)(y) =\frac{1}{2}\langle Ay,y \rangle -\langle p_0,y\rangle +c_0$ with constant  $c_0$,  where $A$ is a linear, symmetric, positive definite operator and $p_0$ is a vector. Noting $\nabla(g_1+f)(y)=Ay-p_0$, the update in \eqref{update_x_bar} can be reformulated as the following preconditioned iteration 
\begin{equation}
    \bar x^n =  y^n + \mathbb{M}^{-1}(p^n - Ay^n),
\end{equation}
where 
\[
 \mathbb{M}= A +M, \quad p^n := p_0 + \xi^n, \ \text{with}\ \xi ^n\in \partial g_2\left( x^n \right)
\]
which is one time classical preconditioned iteration for solving \eqref{eq:topre}, i.e., $Ay = p^n$ with initial value $y^n$.
\end{proposition}
\begin{proof}
    With the assumptions, the equation \eqref{eq:topre} becomes
    \[
    - \xi ^n + (M+A)y - My^n -p_0=0,
    \]
    which is equivalent to 
    \begin{align*}
           y &= (M+A)^{-1}(My^n+\xi^n+p_0) 
         = (M+A)^{-1}((M+A)y^n + \xi^n+p_0 - Ay^n) \\
        & = y^n + (M+A)^{-1}( \xi^n+p_0 - Ay^n)= y^n + (M+A)^{-1}(p^n - Ay^n).
    \end{align*}
    By letting $\bar x^n=y$, the proof is finished. We emphasize that the preconditioned iteration automatically determines $M$, and no explicit expression is needed. We refer to \cite[Section 2.3]{BSCC} for more details and \cite[Proposition 3.2]{Shensun2023} for a similar idea. \qed
\end{proof}

\section{Convergence analysis}\label{sec:conver}
This section analyzes the convergence of the two proposed algorithms and addresses a generalized special case. We first present the convergence for the case that $f(x)$ and $g_1(x)$ have Lipschitz gradients while $g_2$ is non-smooth. For the case $g_1$ being non-smooth, we rigorously prove convergence for the scenario where $g_1(x)$ is non-differentiable, and $g_2(x)$ is Lipschitz differentiable. These two cases can cover the most practical computational scenarios encountered in applications.


\subsection{Convergence analysis of Algorithm \ref{algorithm1}}

We now begin our convergence analysis by establishing the following lemma that illustrates the quasi-descent property inherent in the underlying energy function.

\begin{lemma}\label{lemma1}
Let $\left\{ \bar{x}^n \right\}_n $ be generated by \eqref{eq:algorithmic:update1}, then the following inequality holds
\begin{align}\label{basic_ineq1}
&E\left( \bar{x}^n \right) \le E\left( x^n \right) +\frac{\beta_n^2}{2}\left\| x^n-x^{n-1} \right\| _{M}^{2}-\frac{1}{2}\left\| \bar{x}^n-x^n \right\| _{M}^{2}.
\end{align}      
\end{lemma} 

\begin{proof}
As $\bar{x}^n$ is a minimizer of the strongly convex function  in \eqref{eq:algorithmic:update1}, from the first order optimality condition, we obtain
\begin{equation}\label{eq:optimal:st}
    0 \in -\xi^n + M(\bar{x}^n - y^n) + \partial f(\bar{x}^n) + \partial g_1(\bar{x}^n).
\end{equation}
By convexity of $f$ and $g_1$, for any $x^n$, we have
\begin{align*}
    &f(\bar{x}^n) + \langle s_f, x^n - \bar{x}^n \rangle \leq f(x^n), \quad \forall s_f \in \partial f(\bar{x}^n), \\
    &g_1(\bar{x}^n) + \langle s_{g_1}, x^n - \bar{x}^n \rangle \leq g_1(x^n), \quad \forall s_{g_1} \in \partial g_1(\bar{x}^n),
\end{align*}
which leads to the following inequality by noting $\xi^n - M(\bar{x}^n - y^n) \in \partial f(\bar{x}^n) + \partial g_1(\bar{x}^n)$ from \eqref{eq:optimal:st}
\[
\langle \xi^n - M(\bar{x}^n - y^n), x^n - \bar{x}^n \rangle +f(\bar x^n)+g_1(\bar x^n) \leq f(x^n)+g_1(x^n).
\]
By the expansion of $\frac{1}{2}\|x^n - y^n\|_M^2=\frac{1}{2}\|(x^n -\bar x^n) +(\bar x^n - y^n)\|_M^2$, we obtain
\begin{equation*}
     \frac{1}{2}\|\bar{x}^n - y^n\|_M^2 + \langle M(\bar{x}^n - y^n), x^n - \bar{x}^n \rangle + \frac{1}{2}\|\bar{x}^n - x^n\|_M^2=\frac{1}{2}\|x^n - y^n\|_M^2.
\end{equation*}
Adding the above two inequalities gives exactly 
\begin{align}
&\left< -\xi ^n,\bar{x}^n \right> +\frac{1}{2}\left\| \bar{x}^n-y^n \right\| _{M}^{2}+f\left( \bar{x}^n \right) +g_1\left( \bar{x}^n \right) \le \left< -\xi 
^n,x^n \right> +\frac{1}{2}\left\| x^n-y^n \right\| _{M}^{2} \notag\\
&+f\left( x^n \right) +g_1\left( x^n \right) -\frac{1}{2}\left\| \bar{x}^n-x^n \right\| _{M}^{2}. \label{eq:ener:dec:1}
\end{align}
Given $\xi ^n\in \partial g_2\left( x^n \right) $, the convexity of $g_2$ implies
\begin{equation*}
 g_2\left( \bar{x}^n \right) \geq g_2\left( x^n \right) +\left< \xi ^n,\bar{x}^n-x^n \right>  	\Rightarrow  \langle -\xi^n, \bar x^n -x^n\rangle \geq g_2(x^n)-g_2(\bar x^n).    
\end{equation*}
Substituting $\langle -\xi^n, \bar x^n -x^n\rangle $ with $g_2(x^n)-g_2(\bar x^n)$ into the left-hand side of \eqref{eq:ener:dec:1} yields 
\begin{align*}
&g_2\left( x^n \right) -g_2\left( \bar{x}^n \right) +\frac{1}{2}\left\| \bar{x}^n-y^n \right\| _{M}^{2}+f\left( \bar{x}^n \right) +g_1\left( \bar{x}^n \right) \le \frac{1}{2}\left\| x^n-y^n \right\| _{M}^{2}\\
&+f\left( x^n \right) +g_1\left( x^n \right) -\frac{1}{2}\left\| \bar{x}^n-x^n \right\| _{M}^{2}.
\end{align*}
Rearranging the above terms, we obtain
\begin{align*}
&f\left( \bar{x}^n \right) +g_1\left( \bar{x}^n \right) -g_2\left( \bar{x}^n \right) \le f\left( x^n \right) +g_1\left( x^n \right) -g_2\left( x^n \right) +\frac{1}{2}\left\| x^n-y^n \right\| _{M}^{2} \\
&-\frac{1}{2}\left\| \bar{x}^n-x^n \right\| _{M}^{2}.  
\end{align*}
Remembering $y^n = x^n +\beta_n(x^n-x^{n-1})$, this leads to the inequality  \eqref{basic_ineq1}. \qed
\end{proof}
This important descent characteristic rigorously justifies the adoption of a non-monotone line search strategy. Below, we provide a concise characterization of this mechanism's mathematical foundation, verifying its well-posedness through the following systematic analysis.




\begin{remark}
The relationship between non-monotone and monotone line search in Algorithm \ref{algorithm1} is governed by the parameter $\nu _n$. When $\nu _n=0$ for all iterations $n$, the non-monotone line search reduces exactly to its monotone counterpart, requiring strict energy descent. However, the non-monotone line search also requires the critical descent condition 
\begin{equation*}
\nu_n<\eta\cdot\lambda_n\cdot\left\|d^n\right\|^2, \ \text{for sufficiently large} \ n.    
\end{equation*}
This inequality is readily satisfiable under our parameterization of $\nu _n$. Specifically, for $n>n_0$ (where $n_0$ is problem-dependent), we  obtain
\begin{equation*}
\frac{w}{n+1}<\eta \cdot \lambda _{\max}\cdot \rho^{N_{\max}}.
\end{equation*}
Here $\lambda_{\text{max}}$ denotes the maximum admissible step size and $\rho^{N_\text{max}}$ reflects the worst-case step size contraction. While the line search permits temporary function value increases during initial iterations (enabled by $\nu_n$), the descent dominance condition $\nu_n<\eta\cdot\lambda_n\cdot\left\|d^n\right\|^2$, guarantees eventual monotone convergence. 
Furthermore, this relaxation enables temporary increases in the function value while maintaining global convergence, which is a critical feature for escaping shallow local minima in non-convex landscapes \cite[Example 3.4]{ArtachoSIAM}.
\end{remark}


We now establish global convergence properties for the proposed Algorithm \ref{algorithm1}. This result will serve as the foundation for deriving global convergence of the sequence $x^n$ under appropriate regularity conditions.

\begin{lemma}\label{lemma2}
Let $\left\{ x^n \right\}_n $ is generated from $npDCAe_{nls}$ as Algorithm \ref{algorithm1} for solving \eqref{min_problem}, then we have
\begin{equation}\label{convergence1}
\sum_{n=1}^{\infty}{\left\| x^{n+1}-x^n \right\| ^2<\infty}. 
\end{equation}
\end{lemma}
\begin{proof}
Building upon the inequality \eqref{basic_ineq1}, we now incorporate it with the non-monotone line search. When a suitable step size $\lambda _n >0$ in Case 1 of \eqref{LSDE} satisfies the search criterion, we obtain
\begin{align*}
E\left( x^{n+1} \right) \triangleq E\left( \bar{x}^n+\lambda _nd^n \right) \le E\left( \bar{x}^n \right) -\eta \cdot \lambda _n\cdot \left\| d^n \right\| ^2+\nu _n.    
\end{align*}
Combining these results yields
\begin{align*}
E\left( x^{n+1} \right) \le E\left( x^n \right) -\eta \cdot \lambda _n\cdot \left\| d^n \right\| ^2-\frac{1}{2}\left\| \bar{x}^n-x^n \right\| _{M}^{2}+\frac{\beta _{n}^{2}}{2}\left\| x^n-x^{n-1} \right\| _{M}^{2}+\nu _n. 
\end{align*}
Rearranging terms produces the key inequality
\begin{align}
&\left( \eta \cdot \lambda _n-\frac{\omega}{n+1} \right) \left\| \bar{x}^n-x^n \right\| ^2+\frac{1}{2}\left\| \bar{x}^n-x^n \right\| _{M}^{2}-\frac{\beta _{n}^{2}}{2}\left\| x^n-x^{n-1} \right\| _{M}^{2} \notag \\
&\le E\left( x^n \right) -E\left( x^{n+1} \right) =f\left( x^n \right) +g\left( x^n \right) -f\left( x^{n+1} \right) -g\left( x^{n+1} \right). \label{eq:mix1} 
\end{align}
For Case 2 of \eqref{LSDE}, while the line search does not succeed with $N_{\max}$ searches, we assume $\lambda_n=0$ along with $\nu_n$ set to 0 and hence $x^{n+1} = \bar x^{n}$ as in \eqref{basic_ineq1}, i.e., 
\begin{equation}\label{eq:mix2}
    E\left( {x}^{n+1} \right) \le E\left( x^n \right) +\frac{\beta_n^2}{2}\left\| x^n-x^{n-1} \right\| _{M}^{2}-\frac{1}{2}\left\| {x}^{n+1}-x^n \right\| _{M}^{2}.
\end{equation}
For unifying the analysis of both cases, we first define an indicator function 
\[
\delta _{\lambda}^{n}=\begin{cases}	1   \quad \lambda _n>0\\	0  \quad \lambda _n=0\\\end{cases}. 
\]
We now combine the inequality \eqref{eq:mix1} and \eqref{eq:mix2} with $\delta _{\lambda}^{n}$ as follows
\begin{align*}
&\delta _{\lambda}^{n}\left[ \left( \eta \lambda _n-\frac{\omega}{n+1} \right) \left\| \bar{x}^n-x^n \right\| ^2+\frac{1}{2}\left\| \bar{x}^n-x^n \right\| _{M}^{2}-\frac{\beta _{n}^{2}}{2}\left\| x^n-x^{n-1} \right\| _{M}^{2} \right] \\
&+\left( 1-\delta _{\lambda}^{n} \right) \left[ \frac{1}{2}\left\| x^{n+1}-x^n \right\| _{M}^{2}-\frac{\beta _{n}^{2}}{2}\left\| x^n-x^{n-1} \right\| _{M}^{2} \right] \\
&\le f\left( x^n \right) +g\left( x^n \right) -f\left( x^{n+1} \right) -g\left( x^{n+1} \right). 
\end{align*}
Whether $\lambda _n= 0$ or $\lambda _n\ne 0$, and with the condition $\frac{1}{\left( 1+\lambda _{n-1} \right) ^2}>\beta _{n}^{2}$, we can obtain
\begin{align*}
&\delta _{\lambda}^{n}\left( \eta \lambda _n-\frac{\omega}{n+1} \right) \left\| \bar{x}^n-x^n \right\| ^2+\frac{1}{2}\left( \frac{1}{\left( 1+\lambda _{n-1} \right) ^2}-\beta _{n}^{2} \right) \left\| x^n-x^{n-1} \right\| _{M}^{^2}\\
&-\frac{1}{2\left( 1+\lambda _{n-1} \right) ^2}\left\| x^n-x^{n-1} \right\| _{M}^{^2}+\frac{\delta_{\lambda}^{n}}{2}\left\| \bar{x}^n-x^n \right\| _{M}^{2}+\frac{1-\delta_{\lambda}^{n}}{2}\left\| x^{n+1}-x^n \right\| _{M}^{2}\\
&= \delta _{\lambda}^{n}\left( \eta \lambda _n-\frac{\omega}{n+1} \right) \left\| \bar{x}^n-x^n \right\| ^2-\frac{\beta _{n}^{2}}{2}\left\| x^n-x^{n-1} \right\| _{M}^{2}+\frac{\delta _{\lambda}^{n}}{2}\left\| \bar{x}^n-x^n \right\| _{M}^{2}\\
&+\frac{\left( 1-\delta _{\lambda}^{n} \right)}{2}\left\| x^{n+1}-x^n \right\| _{M}^{2}\le f\left( x^n \right) +g\left( x^n \right) -f\left( x^{n+1} \right) -g\left( x^{n+1} \right) .
\end{align*}
Recalling that $\bar{x}^{n-1}-x^{n-1}=\frac{x^n-x^{n-1}}{1+\lambda _{n-1}}$ along with $\|\bar{x}^{n-1}-x^{n-1}\|_{M}^2= \frac{1}{(1+\lambda _{n-1})^2} \|x^n-x^{n-1}\|_{M}^2$, and noting that
\begin{equation}\label{eq:up:barx:barxn}
\frac{\delta_{\lambda}^{n}}{2}\left\| \bar{x}^n-x^n \right\| _{M}^{2}+\frac{1-\delta_{\lambda}^{n}}{2}\left\| x^{n+1}-x^n \right\| _{M}^{2} \geq \frac{1}{2}\left\| \bar{x}^n-x^n \right\| _{M}^{2}, 
\end{equation}
then we obtain
\begin{align*}
&\delta _{\lambda}^{n}\left( \eta \lambda _n-\frac{\omega}{n+1} \right) \left\| \bar{x}^n-x^n \right\| ^2+\frac{1}{2}\left( \frac{1}{\left( 1+\lambda _{n-1} \right) ^2}-\beta _{n}^{2} \right) \left\| x^n-x^{n-1} \right\| _{M}^{^2}\\
&\le f\left( x^n \right) +g\left( x^n \right) +\frac{1}{2}\left\| \bar{x}^{n-1}-x^{n-1} \right\| _{M}^{2}-f\left( x^{n+1} \right) -g\left( x^{n+1} \right) -\frac{1}{2}\left\| \bar{x}^n-x^n \right\| _{M}^{2}.
\end{align*}
Reformulating the above statements and from Proposition \ref{prop1}, we arrive at 
\begin{align}
&\frac{C_{\lambda}\mu _{\min}\left( M \right)}{2}\left\| x^n-x^{n-1} \right\| ^2\le \frac{C_{\lambda}}{2}\left\| x^n-x^{n-1} \right\| _{M}^{2}\le f\left( x^n \right) +g\left( x^n \right) \\
&+\frac{1}{2}\left\| \bar{x}^{n-1}-x^{n-1} \right\| _{M}^{2}-f\left( x^{n+1} \right) -g\left( x^{n+1} \right) -\frac{1}{2}\left\| \bar{x}^n-x^n \right\| _{M}^{2}.
\label{eq:mix:final}
\end{align}
Summing \eqref{eq:mix:final} from $n_0$ to $N$, we further get
\begin{align*}
&\frac{C_{\lambda}\mu _{\min}\left( M \right)}{2}\cdot \sum_{n=n_0}^N{\left\| x^n-x^{n-1} \right\| ^2}\le f\left( x^{n_0} \right) +g\left( x^{n_0} \right) \\
&+\frac{1}{2}\left\| \bar{x}^{n_0-1}-x^{n_0-1} \right\| _{M}^{2}-f\left( x^{N+1} \right) -g\left( x^{N+1} \right) -\frac{1}{2}\left\| \bar{x}^N-x^N \right\| _{M}^{2}.
\end{align*}
In addition, notice that $f(x)+g(x)$ is proper and bounded below, therefore, $f\left( x^n \right) +g\left( x^n \right) +\frac{1}{2}\left\| \bar{x}^{n-1}-x^{n-1} \right\| _{M}^{2}$ is non-increasing and bounded below. Consequently, letting $N\rightarrow \infty $, we have the following inequality
\begin{align*}
&\frac{C_{\lambda}\mu _{\min}\left( M \right)}{2}\cdot \sum_{n=n_0}^{\infty}{\left\| x^n-x^{n-1} \right\| ^2}\le f\left( x^{n_0} \right) +g\left( x^{n_0} \right) +\frac{1}{2}\left\| \bar{x}^{n_0-1}-x^{n_0-1} \right\| _{M}^{2}\\
&-\mathop {\mathrm{inf}} \limits_{N\rightarrow \infty}\left( f\left( x^{N+1} \right) +g\left( x^{N+1} \right) +\frac{1}{2}\left\| \bar{x}^N-x^N \right\| _{M}^{2} \right) <\infty.
\end{align*}
This leads to \eqref{convergence1}, and the proof is finished. \qed
\end{proof}

We now establish the global convergence of the sequence $x^n$ generated by Algorithm \ref{algorithm1}. This analysis employs a perturbation energy function
\begin{equation} \label{A(x,y)}
A\left(x,\bar{y},y \right) =f\left( x \right) +g\left( x \right) +\frac{1}{2} \left\|\bar{y}-y \right\| _{M}^{2}
\end{equation}
to derive boundedness and convergence guarantees for the iterates. 

\begin{theorem}\label{thm:globalconvergence}
Let $\left\{ x^n \right\}_n $ is generated from $npDCAe_{nls}$ in Algorithm \ref{algorithm1} for solving \eqref{min_problem}, assume that $\nabla g_1(x)$ is Lipschitz continuous with constant $L_{g_1}$, $M$ is a positive definite and linear operator,  and $A\left(x,\bar{y},y \right)$ has the KL property, then the following properties hold
\begin{itemize}
\item[{\rm (i)}]\textbf{Monotonicity}: The sequence $\left\{ A\left( x^{n+1},\bar{x}^{n},x^n \right) \right\}_n $ is monotonically decreasing and converges to some limit: $\lim_{n\rightarrow \infty} A\left( x^{n+1},\bar{x}^{n},x^n \right) =\zeta $. 
\item[{\rm (ii)}]\textbf{Boundedness}\label{Boundedness}: The iterates $\left\{ x^n \right\}_n $ is bounded. 
\item[{\rm (iii)}]\textbf{Optimality}: Any cluster point of $\left\{ x^n \right\}_n $ is the critical point of the problem \eqref{min_problem}, i.e. it satisfies the first order optimality condition. 
\item[{\rm (iv)}]\textbf{Convergence}: The sequence $\left\{ x^n \right\}_n $ is globally convergent with $\sum_{n=1}^{\infty}{\left\| x^{n}-x^{n-1} \right\|}<\infty $.
\end{itemize}
\end{theorem}
\begin{proof}   
(i). 
Setting $C_1=\frac{C_{\lambda}\mu _{\min}\left( M \right)}{2}$, with \eqref{eq:mix:final},  we obtain
\begin{equation}\label{ineq1_Th1}
A\left( x^{n+1},\bar{x}^{n},x^n \right) \le A\left( x^n,\bar{x}^{n-1} x^{n-1} \right)- C_1\left\| x^{n}-x^{n-1} \right\|^2 . \end{equation}
This indicates that $\left\{ A\left( x^{n+1},\bar{x}^{n}, x^n \right) \right\} _{n>n_0}$ is monotonically decreasing, given the lower boundedness of $f\left(x\right)+g\left(x\right)$. Consequently,  we have
\begin{equation}\label{eq1_Th1}
\lim_{n\rightarrow \infty} A\left( x^{n+1},\bar{x}^{n}, x^n \right) =\zeta.   
\end{equation}

(ii). We know that the monotone decrease of $\left\{ A\left( x^{n+1},\bar{x}^{n}, x^n \right)\right\}_{n>n_0}$ established in (i) implies
\begin{equation*}
f\left( x^{n+1} \right) +g\left( x^{n+1} \right) \le A\left( x^{n+1},\bar{x}^n,x^n \right) \le f\left( x^{n_0+1} \right) +g\left( x^{n_0+1} \right) +\frac{1}{2}\left\| \bar{x}^{n_0}-x^{n_0} \right\| _{M}^{2}, \ \forall n > n_0
\end{equation*}
The level-boundedness of $f\left(x\right)+g\left(x\right)$ together with \eqref{convergence1} guarantees the boundedness of $\left\{x^n\right\}_{n}$. 

(iii). Given the boundedness of $\left\{x^n\right\}_n$ from (ii), there exists at least one cluster point $x^*$ with a convergent subsequence $\left\{x^{n_i}\right\}_{n_i}$ satisfying
\begin{equation*}
\lim_{i\rightarrow \infty} x^{n_i}=x^*\text{ and } \lim_{i\rightarrow \infty}\left\| x^{n_i}-x^{n_i-1} \right\|\rightarrow0.
\end{equation*}
The momentum update $y^{n_i}=\left( x^{n_i}+\beta_{n_i}\left( x^{n_i}-x^{n_i-1} \right) \right)$ hence satisfies
\begin{equation}
\lim_{i\rightarrow \infty} y^{n_i}=x^*
\end{equation}
due to the boundedness of $\left\{\beta _n\right\}_n$ and vanishing step differences \eqref{convergence1}.
From the first-order optimality condition of the proximal update \eqref{eq:topre} at $y = \bar x^{n_i}$, we have
\[
M\left(\bar{x}^{n_i}-y^{n_i}\right) +\nabla f\left(\bar{x}^{n_i}\right) +\nabla g_1\left(\bar{x}^{n_i}\right) \in \partial g_{2}({x}^{n_i}).
\]
We thus arrive at
\[
 -M\left(\bar{x}^{n_i}-y^{n_i}\right) -\nabla f\left(\bar{x}^{n_i}\right) -\nabla g_1\left(\bar{x}^{n_i}\right)  +\nabla (f+g_1)({x}^{n_i})  \in -\partial g_2({x}^{n_i}) +\nabla (f+g_1)({x}^{n_i}). 
\]
The left-hand side of the above equation can be written as
\begin{equation}\label{eq:lefth1}
-M\left(\bar{x}^{n_i}-y^{n_i}\right) + (\nabla f({x}^{n_i})- \nabla f(\bar {x}^{n_i})) + (\nabla g_1({x}^{n_i})- \nabla g_1(\bar {x}^{n_i})). 
\end{equation}
With $\bar x^{n_i}-x^{n_i} = (x^{n_i +1}-x^{n_i})/(1+\lambda_{n_i})$ and 
$\bar{x}^{n_i}-y^{n_i}=\bar x^{n_i}-x^{n_i} - \beta_{n_i}(x^{n_i}-x^{n_i -1}) $, the Lipschitz continuity of $\nabla f$ and $\nabla g_1$, and \eqref{convergence1}, we conclude that the formulation in \eqref{eq:lefth1} tends to zero while $n_i \rightarrow \infty$. We then obtain
\[
0 \in \lim_{n_i \rightarrow \infty } -\partial g_2({x}^{n_i}) +\nabla (f+g_1)({x}^{n_i}).
\]
Remembering $x^{n_i} \rightarrow x^*$, with closedness of $\partial g_2$ \cite[Proposition 16.36]{HBPL} or \cite[Theorem 8.6]{Roc1}, we have 
\begin{align*}
0\in -\partial g_2\left( x^* \right) +\nabla f\left( x^* \right) +\nabla  g_1\left( x^* \right),    
\end{align*}
confirming that $x^*$ satisfies the first-order optimality condition for \eqref{min_problem}. 

(iv). The subgradient of the energy function is given by
\begin{equation}
\partial A\left( x,\bar{y},y \right) =\left( \begin{array}{c}	\nabla f\left( x \right) +\nabla g_1\left( x \right)-\partial g_2(x)\\
M\left( \bar{y}-y\right)\\	-M\left( \bar{y}-y\right)\\\end{array} \right).
\end{equation}
From the optimality condition
\begin{equation*}
0\in -\xi ^n+M\left( \bar{x}^n-y^n \right) +\nabla f\left( \bar{x}^n \right) +\nabla g_1\left( \bar{x}^n \right),
\end{equation*}
we derive the inequality
\begin{align*}
&\mathrm{dist}\left(0,\partial A\left( x^n,\bar{x}^{n-1},x^{n-1} \right) \right) \le \left( L+L_{g_1} \right) \left\| \bar{x}^n-x^n \right\| +\left\| M\left( \bar{x}^n-y^n \right) \right\| +\left\|M(\bar{x}^{n-1}-x^{n-1}) \right\|.
\end{align*} 
Setting $\tilde{L}=L+L_{g_1}$, according to $\bar{x}^n-x^n=\frac{\lambda _n}{1+\lambda _n}\left( x^{n+1}-x^n \right) $, $\bar{x}^{n-1}-x^{n-1}=\frac{\lambda_{n-1}}{1+\lambda_{n-1}}\left(x^n-x^{n-1} \right)$, 
we obtain the following formulation
\begin{align*}
&\mathrm{dist}\left(0,\partial A\left( x^n,\bar{x}^{n-1},x^{n-1} \right) \right) \le  \tilde{L}\left\| x^{n+1}-x^n \right\| +\left\| M\left( x^{n+1}-x^n \right) \right\| \\
&+\beta _n\left\| M\left( x^n-x^{n-1} \right) \right\| +\left\|M\left( x^n-x^{n-1} \right) \right\|.
\end{align*}
Together with the boundedness of $M$ and $\beta _n$, we have 
\begin{align*}
\mathrm{dist}\left(0,\partial A\left( x^n,\bar{x}^{n-1},x^{n-1} \right) \right) \le C\left( \left\| x^{n+1}-x^n \right\| +\left\| x^n-x^{n-1} \right\| \right),
\end{align*}
where $C$ is a constant. With $\lim_{n\rightarrow \infty} \left\| x^n-x^{n-1} \right\| \rightarrow 0$ by \eqref{convergence1},  we obtain
\begin{equation}
\lim_{n\rightarrow \infty} \mathrm{dist}\left(0,\partial A\left( x^n,\bar{x}^{n-1},x^{n-1} \right) \right) = 0.    
\end{equation}
Given the KL property of $A\left(x,\bar{y},y \right)$ with concave function $\psi\left(\cdot\right)$, and from \eqref{eq1_Th1}, we establish
\begin{equation*}
\begin{aligned}
&\left[ \psi \left( A\left( x^n,\bar{x}^{n-1},x^{n-1} \right) -\zeta \right) -\psi \left( A\left( x^{n+1},\bar{x}^{n},x^n \right) -\zeta \right) \right] \cdot \mathrm{dist}\left(0, \partial A\left( x^n,\bar{x}^{n-1},x^{n-1} \right) \right) \\
&\ge \psi \prime\left( A\left( x^n,\bar{x}^{n-1},x^{n-1} \right) -\zeta \right) \cdot \left[ A\left( x^n,\bar{x}^{n-1},x^{n-1} \right)-A\left( x^{n+1},\bar{x}^{n},x^n \right)  \right] \cdot \\
&\mathrm{dist}\left(0,\partial A\left( x^n,\bar{x}^{n-1},x^{n-1} \right) \right) \ge A\left( x^n,\bar{x}^{n-1},x^{n-1} \right) -A\left( x^{n+1},\bar{x}^{n},x^n \right) \ge C_1\left\| x^{n}-x^{n-1} \right\| ^2.
\end{aligned}   
\end{equation*}
Set $\varPsi \left( x^{n-1},\bar{x}^{n-1}, x^n,\bar{x}^{n},x^{n+1},\zeta \right) =\psi \left( A\left( x^n,\bar{x}^{n-1},x^{n-1} \right) -\zeta \right) -\psi \left( A\left( x^{n+1},\bar{x}^{n},x^n \right) -\zeta \right) $.
We further obtain the following statement
\begin{align*}
\left\| x^{n}-x^{n-1} \right\| ^2 \le \tilde{C}\cdot \varPsi \left( x^{n-1},\bar{x}^{n-1}x^n,\bar{x}^{n},x^{n+1},\zeta \right) \cdot \left( \left\| x^{n+1}-x^n \right\| +\left\| x^n-x^{n-1} \right\| \right),  
\end{align*}
where $\tilde{C}=\frac{C}{C_1}$. From the basic inequality $a^2\le cd\Longrightarrow a\le c+\frac{d}{4}$,  we obtain
\begin{align}\label{ineq:splitform}
&\frac{1}{2}\left\| x^{n}-x^{n-1} \right\| \le \tilde{C}\cdot \left[ \psi \left( A\left( x^n,\bar{x}^{n-1},x^{n-1} \right) -\zeta \right) -\psi \left( A\left( x^{n+1},\bar{x}^{n},x^n \right) -\zeta \right) \right] \\
&+\frac{1}{4}\left( \left\| x^{n+1}-x^n \right\| -\left\| x^{n}-x^{n-1} \right\| \right). 
\end{align}
Summing both sides of the inequality, we obtain
\begin{align*}
&\frac{1}{2}\sum_{n=n_0}^N{\left\| x^n-x^{n-1} \right\| \le \tilde{C}}\left[ \psi \left( A\left( x^{n_0},\bar{x}^{n_0-1},x^{n_0-1} \right) -\zeta \right) -\psi \left( A\left( x^{N+1},\bar{x}^{N},x^N \right) -\zeta \right) \right] \\
&+\frac{1}{4}\left( \left\| x^{N+1}-x^{N} \right\| -\left\| x^{n_0}-x^{n_0-1} \right\| \right). 
\end{align*}
Letting $N\rightarrow \infty $,  we know that  
\begin{equation*}
\lim_{N\rightarrow \infty} \psi \left( A\left( x^{N+1},\bar{x}^{N},x^N \right)  -\zeta \right) =\psi \left( 0 \right) =0.    
\end{equation*}
Therefore, with \eqref{convergence1}, we obtain the desired statement
\begin{align*}
&\sum_{n=n_0}^{\infty}{\left\| x^{n}-x^{n-1} \right\| \le 2\tilde{C}\cdot \psi \left( A\left( x^{n_0},\bar{x}^{n_0-1},x^{n_0-1} \right) -\zeta \right)}< \infty.  
\end{align*}
Finally, we obtain
$\sum_{n=1}^{\infty}{\left\| x^{n}-x^{n-1} \right\|}< \infty$. \qed
\end{proof}

We now establish convergence properties for Algorithm \ref{algorithm2}. While Algorithm \ref{algorithm2} differs from Algorithm \ref{algorithm1} only in its DC decomposition strategy, their convergence proofs share identical methodology. We present three parallel conclusions for Algorithm \ref{algorithm2} for completeness.

\subsection{Convergence analysis of Algorithm \ref{algorithm2}}
The following Lemma \ref{lemma3}, Lemma \ref{lemma4}, and Theorem \ref{th2} are obtained on this basis, and we will not write them out one by one later on. Following the analytical framework established for Algorithm \ref{algorithm2}, we initially demonstrate that Lemma \ref{lemma3} serves to rigorously characterize the quasi-descent property intrinsic to the energy functional governing the system.
\begin{lemma} \label{lemma3}
Let $\{\bar{x}^n\}_n$ is generated by \eqref{eq:algorithmic:update2}, then the following inequality holds
\begin{equation}\label{basic_ineq2}
E\left( \bar{x}^n \right) \le E\left( x^n \right) +\frac{L\beta_n^2}{2}\left\| x^n-x^{n-1} \right\| _{M}^{2}-\frac{L}{2}\left\| \bar{x}^n-x^n \right\| _{M}^{2}.    
\end{equation}
\end{lemma}
\begin{proof}
Firstly, from \eqref{eq:algorithmic:update2},  we know that $\bar{x}^n$ is the minimum point of the strongly convex function
\begin{equation*}
\left< \nabla f\left( y^n \right) -\xi ^n,y \right> +\frac{L}{2}\left\| y-y^n \right\| _{M}^{2}+g_1\left( y \right).  
\end{equation*}
This implies the optimality inequality
\begin{align*}
&\left< \nabla f\left( y^n \right) -\xi ^n,\bar{x}^n \right> +\frac{L}{2}\left\| \bar{x}^n-y^n \right\| _{M}^{2}+g_1\left( \bar{x}^n \right) \le \left< \nabla f\left( y^n \right) -\xi ^n,x^n \right> +\frac{L}{2}\left\| x^n-y^n \right\| _{M}^{2} \\
&+g_1\left( x^n \right) -\frac{L}{2}\left\| \bar{x}^n-x^n \right\| _{M}^{2}.
\end{align*}
Using the subgradient property $\xi ^n\in \partial g_2\left( x^n \right) $, then we obtain
\begin{equation*}
g_2\left( \bar{x}^n \right) \geq g_2\left( x^n \right) +\left< \xi ^n,\bar{x}^n-x^n \right>.    
\end{equation*}
Combining these results yields
\begin{align}
&\left< \nabla f\left( y^n \right), \bar{x}^n-x^n \right> +\frac{L}{2}\left\| \bar{x}^n-y^n \right\|_{M}^{2}+g_1\left( \bar{x}^n \right) -g_2\left( \bar{x}^n \right) \le \frac{L}{2}\left\| x^n-y^n \right\| _{M}^{2} \notag \\
&+g_1\left( x^n \right) 
-g_2\left( x^n \right) -\frac{L}{2}\left\| \bar{x}^n-x^n \right\| _{M}^{2}.  \notag  
\end{align}
Since $f(x)$ is convex and $\nabla f(x)$ is Lipschitz continuous, we have the following two inequalities
\begin{align}
&f\left( \bar{x}^n \right) \le f\left( y^n \right) +\left< \nabla f \left( y^n \right), \bar{x}^n-y^n \right> + \frac{L}{2}\left\| \bar{x}^n-y^n \right\| ^2, \notag  \\
&f\left( x^n \right) \ge f\left( y^n \right) +\left < \nabla f \left( y^n \right), x^n-y^n \right>. \notag
\end{align}
Sum the above two inequalities, and obtain
\begin{align*}
&f\left( \bar{x}^n \right) -g_2\left( \bar{x}^n \right) +\frac{L}{2}\left\| \bar{x}^n-y^n \right\| _{M}^{2}+g_1\left( \bar{x}^n \right) \le \frac{L}{2}\left\| x^n-y^n \right\| _{M}^{2}+g_1\left( x^n \right). \\
&-\frac{L}{2}\left\| \bar{x}^n-x^n \right\| _{M}^{2}-g_2\left( x^n \right) +f\left( x^n \right) +\frac{L}{2}\left\| \bar{x}^n-y^n \right\| ^2.
\end{align*}
Applying the matrix dominance $M\succeq I$ and $y^n = x^n+\beta_n(x^n-x^{n-1})$,  we derive \eqref{basic_ineq2}.\qed
\end{proof}

In the subsequent analysis, we establish the global subsequential convergence property for Algorithm \ref{algorithm2}. To provide a basis for this study, we first present the following preparatory Lemma \ref{lemma4}, whose proof procedure is essentially similar to that of Algorithm \ref{algorithm1}, and we omit it here.

\begin{lemma}\label{lemma4}
Let $\left\{ x^n \right\}_n $ is generated from $pDCAe_{\text{nls}}$ in Algorithm \ref{algorithm2} for solving \eqref{min_problem}, then we have
\begin{equation}\label{convergence2}
\sum_{n=1}^{\infty}{\left\| x^{n+1}-x^n \right\| ^2<\infty}.
\end{equation}
\end{lemma}

To conclude this subsection, we establish the global convergence of Algorithm \ref{algorithm2}. For consistency, we adopt identical conventions as those defined for Algorithm \ref{algorithm1}, given the essential similarity in their analytical frameworks. The detailed theoretical proof of this lemma parallels the methodology employed for Algorithm \ref{algorithm1}.

\begin{theorem}\label{th2} Let $\left\{ x^n \right\}_n $ is generated from pDCAe$_{\text{nls}}$ in Algorithm \ref{algorithm2} for solving \eqref{min_problem}. Define the perturbation energy function
\begin{equation*}
A\left(x,\bar{y},y\right) =f\left(x\right) +g\left(x\right) +\frac{L}{2}\left\|\bar{y}-y\right\|_{M}^{2}.
\end{equation*}
Assume that $\nabla g_1(x)$ is Lipschitz continuous with constant $L_{g_1}$, $M\succeq I$ being a positive definite and linear operator, and $A\left( x,\bar{y},y \right)$ has the KL property, then the following properties hold
\begin{itemize}
\item[{\rm (i)}]\textbf{Monotonicity}: The sequence $\left\{ A\left( x^{n+1},\bar{x}^{n},x^n \right) \right\}_n $ is monotonically decreasing and converges to some limit: $\lim_{n\rightarrow \infty} A\left( x^{n+1},\bar{x}^{n},x^n \right) =\zeta$. 
\item[{\rm (ii)}]\textbf{Boundedness}: The iterates $\left\{ x^n \right\}_n $ is bounded. 
\item[{\rm (iii)}]\textbf{Optimality}: Any cluster point of $\left\{ x^n \right\}_n $ is the critical point of the problem \eqref{min_problem}, i.e. it satisfies the first order optimality condition. 
\item[{\rm (iv)}]\textbf{Convergence}:  The sequence $\left\{ x^n \right\}_n $ is globally convergent with $\sum_{n=1}^{\infty}{\left\| x^{n}-x^{n-1} \right\|}<\infty $.
\end{itemize}
\end{theorem}

\subsection{Convergence analysis for the special case with nonsmooth $g_1$}\label{sec:nonsmooh}
Throughout this paper, suppose that  $g_1(x) $ is proper, lower semicontinuous and convex as follows
\begin{equation}\label{local_lip2}
    g_1\left( x \right) \ge g_1\left( \bar{x} \right) +\langle x-\bar{x},\tilde{w} \rangle +\frac{\mu}{2}\left\| x-\bar{x} \right\|^2, \quad \mu \geq 0,
\end{equation}
for any $\bar x \in \dom g_1$ and $\tilde  w \in \partial g_1(\bar x)$ and $\mu$ can be zero while $g_1$ is not strongly convex. 
We propose a more general case which can cover the widely used case that $g_1(x)$ is non-smooth, along with $\nabla g_2(x)$ being Lipschitz continuous with constant $L_{g_2}$, which includes some important regularization such as SCAD as discussed in \cite{Apx}. 

Now we introduce the following perturbation energy function
\begin{equation} \label{H(x)}
\widetilde{H}\left( \bar{x},x \right) =E\left( \bar{x} \right) +\frac{1}{2}\left\| \bar{x}-x \right\| _{M}^{2}=f\left( \bar{x} \right) +g\left( \bar{x} \right) +\frac{1}{2}\left\| \bar{x}-x \right\| _{M}^{2} ,
\end{equation}
for proving the global convergence of the Algorithm \ref{algorithm1} for the case $g_1$ being non-smooth. For the KL analysis with $\partial \widetilde{H}( \bar{x}^n,x^n )$, it will turn out that the most challenging problem is caused by $\partial g_1(\bar x^n)$. With $\tilde H(\bar x, x)$, one can replace $\partial g_1(\bar x^n)$ through the optimality condition \eqref{eq:topre}. 
\begin{theorem}\label{thm:3:line}
Let $\left\{ x^n \right\} _n$ is generated from $npDCAe_{nls}$ in Algorithm \ref{algorithm1} for solving \eqref{min_problem}. Supposing that $\left\{ x^n \right\} _n$ has a cluster $x^*$ and $\widetilde{H}\left( \bar{x},x \right) $ is a KL function, $N_0=\max\mathrm{(}n_0,n_1)$, the following statements hold
\begin{itemize}
\item[{\rm (i)}] \textbf{Monotonicity and boundedness}: The sequence $\left\{ \widetilde{H}\left( \bar{x}^n,x^n \right) \right\} _{n\ge N_0}$ is decreasing monotonically and converges to some limits: $\lim_{n\rightarrow \infty} \widetilde{H}\left( \bar{x}^n,x^n \right) =\bar{\zeta}$. Moreover, the iteration sequence $\{x^n\}_n$ is bounded.
\item[{\rm (ii)}] \textbf{Optimality}: Any cluster point of $\left\{ x^n \right\}_n $ is a critical point of the problem \eqref{min_problem}, i.e. it satisfies the first order optimality condition.
\item[{\rm (iii)}]  \textbf{Convergence}:  The sequence $\left\{ x^n \right\}_n $ is globally convergent with $\sum_{n=1}^{\infty}{\left\| x^{n}-x^{n-1} \right\|}<\infty $.
\end{itemize}
\end{theorem}

\begin{proof}
(i). From Lemma \ref{lemma1} we derive
\begin{equation*}
E\left( \bar{x}^n \right) \le E\left( x^n \right) +\frac{\beta _{n}^{2}}{2}\left\| x^n-x^{n-1} \right\| _{M}^{2}-\frac{1}{2}\left\| \bar{x}^n-x^n \right\| _{M}^{2}.
\end{equation*}
Due to the non-monotone line search, we have
\begin{equation*}
E\left( x^n \right) \triangleq E\left( \bar{x}^{n-1}+\lambda _{n-1}d^{n-1} \right) \le E\left( \bar{x}^{n-1} \right) -\eta \lambda _{n-1}\left\| d^{n-1} \right\| ^2+\nu _{n-1}.
\end{equation*}
Recalling the definition of $\widetilde{H}\left( \bar{x}^n,x^n \right)$, then by $d^{n-1}=\bar{x}^{n-1}-x^{n-1}$, $\nu _{n-1}=\frac{\omega}{n}\left\| d^{n-1} \right\| ^2$, and  $\bar{x}^{n-1}-x^{n-1}=\frac{x^n-x^{n-1}}{1+\lambda _{n-1}}$, we have
\begin{align*}
&\widetilde{H}\left( \bar{x}^n,x^n \right) =E\left( \bar{x}^n \right) +\frac{1}{2}\left\| \bar{x}^n-x^n \right\| _{M}^{2}\le E\left( x^n \right)+\frac{\beta _{n}^{2}}{2}\left\| x^n-x^{n-1} \right\| _{M}^{2}\\
&\le E\left( \bar{x}^{n-1} \right) -\eta \lambda _{n-1}\left\| d^{n-1} \right\| ^2+\nu _{n-1}+\frac{\beta _{n}^{2}}{2}\left\| x^n-x^{n-1} \right\| _{M}^{2}.
\end{align*}
Noticing that $\widetilde{H}\left( \bar{x}^{n-1},x^{n-1} \right) =E\left( \bar{x}^{n-1} \right) +\frac{1}{2}\left\| \bar{x}^{n-1}-x^{n-1} \right\| _{M}^{2}$, we thus obtain 
\begin{align*}
&\widetilde{H}\left( \bar{x}^{n-1},x^{n-1} \right) -\widetilde{H}\left( \bar{x}^n,x^n \right) \ge \frac{1}{2}\left\| \bar{x}^{n-1}-x^{n-1} \right\| _{M}^{2}+\eta \lambda _{n-1}\left\| d^{n-1} \right\| ^2-\nu _{n-1}-\frac{\beta _{n}^{2}}{2}\left\| x^n-x^{n-1} \right\| _{M}^{2}\\
&\ge \frac{1}{2\left( 1+\lambda _{n-1} \right) ^2}\left\| x^n-x^{n-1} \right\| _{M}^{2}+\frac{\left( \eta \lambda _{n-1}-\frac{\omega}{n} \right)}{\left( 1+\lambda _{n-1} \right) ^2}\left\| x^n-x^{n-1} \right\| ^2-\frac{\beta _{n}^{2}}{2}\left\| x^n-x^{n-1} \right\| _{M}^{2}.
\end{align*}
Since $\eta \lambda _{n-1}-\frac{\omega}{n}>0$, and from the LSDE updates in Proposition \ref{prop1}, we have $\frac{1}{\left( 1+\lambda _{n-1} \right) ^2}-\beta _{n}^{2}>C_{\lambda}$. 
We thus obtain
\begin{equation}
\widetilde{H}\left( \bar{x}^n,x^n \right) -\widetilde{H}\left( \bar{x}^{n-1},x^{n-1} \right) \le -\frac{C_{\lambda}}{2}\left\| x^n-x^{n-1} \right\| _{M}^{2}, 
\end{equation}
which means $\left\{ \widetilde{H}\left( \bar{x}^n,x^n \right) \right\} _{n>N_0}$  is decreasing monotonically. Furthermore, as we see that $f(x)+g(x)$ is bounded below, which implies 
\begin{equation}\label{eq1_Th2}
\lim_{n\rightarrow \infty} \widetilde{H}\left( \bar{x}^n,x^n \right) =\bar{\zeta}.   
\end{equation}
On the other hand, since $f(x)+g(x)$ is proper, closed, and level-bounded,  we conclude that the iteration sequence $\{x^n\}_n$ is bounded with monotone decreasing of $ \widetilde{H}\left(\bar{x}^n,x^n \right)$. 

(ii). Since there exists a cluster point $x^*$ with a convergent subsequence $\left\{x^{n_i}\right\}_{n_i}$ satisfying
\begin{equation*}
\lim_{i\rightarrow \infty} x^{n_i}=x^*\text{ and } \lim_{i\rightarrow \infty}\left\| x^{n_i}-x^{n_i-1} \right\|\rightarrow 0.
\end{equation*}
The momentum update $y^{n_i}=\left( x^{n_i}+\beta ^{n_i}\left( x^{n_i}-x^{n_i-1} \right) \right)$ consequently satisfies 
$\lim_{i\rightarrow \infty} y^{n_i}=x^*$, 
due to the boundedness of $\left\{\beta _n\right\}_n$ and vanishing step differences in \eqref{convergence1}.
From the first-order optimality condition \eqref{eq:topre} at $y = \bar x^{n_i}$, we have
\begin{align*}
0\in -\xi^{n_i}+M\left(\bar{x}^{n_i}-y^{n_i}\right) +\nabla f\left(\bar{x}^{n_i}\right) +\partial g_1\left(\bar{x}^{n_i}\right).
\end{align*}
With $\xi^{n_i} = \nabla g_2(x^{n_i})$, we thus obtain
\[-M\left( \bar{x}^{n_i}-y^{n_i} \right) -\nabla g_2\left( \bar{x}^{n_i} \right) +\nabla g_2\left( x^{n_i} \right) \in \nabla f\left( \bar{x}^{n_i} \right) +\partial g_1\left( \bar{x}^{n_i} \right) -\nabla g_2(\bar{x}^{n_i}).\]
The left-hand side of the above equation can be written as
\begin{equation}\label{eq:lefth2}
-M\left( \bar{x}^{n_i}-y^{n_i} \right) +(\nabla g_2\left( x^{n_i} \right)-\nabla g_2\left( \bar{x}^{n_i} \right)). 
\end{equation}
Furthermore, the cumulative step bound \eqref{convergence1} implies
\begin{equation}\label{eq:xni_cnv}
\lim_{i\rightarrow \infty} \left\|\bar{x}^{n_i}-x^{n_i}\right\|\rightarrow0 \text{ and } \lim_{i\rightarrow \infty} \bar{x}^{n_i}=x^*. 
\end{equation}
With $\bar x^{n_i}-x^{n_i} = (x^{n_i +1}-x^{n_i})/(1+\lambda_{n_i})$ and 
$\bar{x}^{n_i}-y^{n_i}=\bar x^{n_i}-x^{n_i} - \beta_{n_i}(x^{n_i}-x^{n_i -1}) $, the Lipschitz continuity of $\nabla g_2$, and \eqref{eq:xni_cnv}, we conclude that the formulation in \eqref{eq:lefth2} tends to zero while $n_i \rightarrow \infty$. We then obtain the following.
\[0 \in \lim_{n_i \rightarrow \infty } (-\partial g_1(\bar{x}^{n_i}) +\nabla (f+g_2)(\bar{x}^{n_i})).\]
With closedness of $\partial g_1$ \cite[Proposition 16.36]{HBPL} or \cite[Theorem 8.6]{Roc1} and \eqref{eq:xni_cnv}, we have 
\begin{align*}
0\in -\partial g_1\left( x^* \right) +\nabla f\left( x^* \right) +\nabla  g_2\left( x^* \right),    
\end{align*}
confirming $x^*$ satisfies the first-order optimality condition for \eqref{min_problem}. 

(iii). The subgradient of the energy function is given by
\begin{equation}\label{eq:subgrd:H}
\partial \widetilde{H}\left( \bar{x}^n,x^n \right) =\left( \begin{array}{c}	\nabla f\left( \bar{x}^n \right) +\partial g_1\left( \bar{x}^n \right) -\nabla g_2\left( \bar{x}^n \right) +M\left( \bar{x}^n-x^n \right)\\	-M\left( \bar{x}^n-x^n \right)\\\end{array} \right).
\end{equation}
From the first-order optimality condition, we have 
\begin{equation}\label{eq:opti:bar}
0\in \nabla f\left( \bar{x}^{n} \right) +\partial g_1\left( \bar{x}^{n} \right) -\nabla g_2\left( x^{n} \right) +M\left( \bar{x}^{n}-y^{n} \right).
\end{equation}
Substituting
\begin{equation}\label{eq:subs:g1}
    \partial g_1\left( \bar{x}^n \right) \in -\nabla f\left( \bar{x}^n \right) +\nabla g_2\left( x^n \right) -M\left( \bar{x}^n-y^n \right)
\end{equation}
with conditions 
$\bar{x}^{n-1}-x^{n-1}=\frac{x^n-x^{n-1}}{1+\lambda _{n-1}}$ and $x^n-\bar{x}^{n-1}=\lambda _{n-1}\left( \bar{x}^{n-1}-x^{n-1} \right) $, then reformulating the statement, we thus have
\begin{align*}
&\dist\left( 0,\partial \widetilde{H}\left( \bar{x}^n,x^n \right) \right) \le L_{g_2}\left\| \bar{x}^n-x^n \right\| +B_1\left\| \beta _n\left( x^n-x^{n-1} \right) \right\| +B_1\left\| \bar{x}^n-x^n \right\|\\
& \le \left( L_{g_2}+B_1 \right) \left\| x^{n+1}-x^n \right\| +B_1\left\| x^n-x^{n-1} \right\| \le B\left( \left\| x^{n+1}-x^n \right\| +\left\| x^n-x^{n-1} \right\| \right),
\end{align*}
where $B_1$ is the maximum singular value of $M$, $B_2=B_1+L_{g_2}$ and $B=\max(B_1,B_2)$ are constants. From \eqref{convergence1}, we see $\lim_{n\rightarrow \infty} \left\| x^n-x^{n-1} \right\| \rightarrow 0$.  We thus arrive at
\begin{equation*}
\lim_{n\rightarrow \infty} \mathrm{dist}\left( 0,\partial \widetilde{H} \left( \bar{x}^n,x^n \right) \right) =0.
\end{equation*}
Given the KL property of $\widetilde{H} \left( \bar{x},x \right) $ with concave function $\psi \left( \cdot \right) $, and from \eqref{eq1_Th1}, we establish
\begin{align*}
&\left[ \psi \left( \widetilde{H} \left( \bar{x}^{n-1},x^{n-1} \right) -\bar{\zeta}\right) -\psi \left( \widetilde{H} \left( \bar{x}^n,x^n \right) -\bar{\zeta} \right) \right] \cdot \mathrm{dist}\left( 0,\partial \widetilde{H}\left( \bar{x}^{n-1},x^{n-1} \right) \right) \\
&\ge \psi \prime\left( \widetilde{H}\left( \bar{x}^{n-1},x^{n-1} \right) -\bar{\zeta}\right) \cdot \left[\widetilde{H}\left( \bar{x}^{n-1},x^{n-1} \right) - \widetilde{H}\left( \bar{x}^n,x^n \right) \right] \cdot \mathrm{dist}\left( 0,\partial \widetilde{H}\left( \bar{x}^{n-1},x^{n-1} \right) \right) \\
&\ge \widetilde{H}\left( \bar{x}^{n-1},x^{n-1} \right) -\widetilde{H}\left( \bar{x}^n,x^n \right) \ge \frac{C_{\lambda}}{2}\left\| x^n-x^{n-1} \right\| _{M}^{2}\ge C_1\left\| x^n-x^{n-1} \right\| ^2.
\end{align*}
The existence of the positive constant $C_1$ is due to the positive definiteness of $M$. Then we set 
\begin{align*}
\varPsi \left( x^{n-1},\bar{x}^{n-1},x^n,\bar{x}^n,\bar{\zeta} \right) =\psi \left( \widetilde{H}\left( \bar{x}^{n-1},x^{n-1} \right) -\bar{\zeta} \right) -\psi \left( \widetilde{H}\left( \bar{x}^n,x^n \right) -\bar{\zeta} \right).
\end{align*}
We further obtain the following statement
\begin{align*}
\left\| x^n-x^{n-1} \right\| ^2\le \tilde{B}\cdot \varPsi \left( x^{n-1},\bar{x}^{n-1},x^n,\bar{x}^n,\bar{\zeta} \right) \cdot \left( \left\| x^n-x^{n-1} \right\| +\left\| x^{n-1}-x^{n-2} \right\| \right),
\end{align*}
where $\tilde{B}=\frac{B}{C_1}$. From the basic inequality $a^2\le cd\Longrightarrow a\le c+\frac{d}{4}$,  we obtain
\begin{align*}
&\frac{1}{2}\left\| x^n-x^{n-1} \right\| \le \tilde{B}\cdot \left[ \psi \left( \widetilde{H}\left( \bar{x}^{n-1},x^{n-1} \right) -\bar{\zeta} \right) -\psi \left( \widetilde{H}\left( \bar{x}^n,x^n \right) -\bar{\zeta} \right) \right] \\
&+\frac{1}{4}\left( \left\| x^{n-1}-x^{n-2} \right\| -\left\| x^n-x^{n-1} \right\| \right) .
\end{align*}
Summing both sides of the inequality, we obtain
\begin{align*}
&\frac{1}{2}\sum_{n=N_0}^N{\left\| x^n-x^{n-1} \right\|}\le \tilde{B}\left[ \psi \left( \widetilde{H}\left( \bar{x}^{N_0-1},x^{N_0-1} \right) -\bar{\zeta} \right) -\psi \left( \widetilde{H}\left( \bar{x}^N,x^N \right) -\bar{\zeta} \right) \right] \\
&+\frac{1}{4}\left( \left\| x^{N_0-1}-x^{N_0-2} \right\| -\left\| x^N-x^{N-1} \right\| \right) .
\end{align*}
Letting $N\rightarrow \infty $, we know that  
\begin{equation*}
\lim_{N\rightarrow \infty} \psi \left( \widetilde{H}\left( \bar{x}^N,x^N \right) -\bar{\zeta} \right) =\psi \left( 0 \right) =0. 
\end{equation*}
Therefore, from \eqref{convergence1}, we obtain the desired statement 
\begin{align*}
\sum_{n=N_0}^{\infty}{\parallel x^n}-x^{n-1}\parallel \le 2\tilde{B}\cdot \psi  \left( \widetilde{H} \left( \bar{x}^{N_0-1}x^{N_0-1} \right) -\bar{\zeta}\right)+\frac{1}{2}\left\| x^{N_0-1}-x^{N_0-2} \right\| <\infty .   
\end{align*}
Consequently, we obtain
$\sum_{n=1}^{\infty}{\left\| x^{n}-x^{n-1} \right\|}< \infty$ and the global convergence of $\{x^n\}_n$. 
\qed
\end{proof}
Below we give the conclusion of this particular calculus in the framework of Algorithm \ref{algorithm2} and omit a similar proof for compactness.
\begin{theorem}\label{thm:line:ano}
Assume the sequence $\left\{ x^n \right\}_n $ is generated from $pDCAe_{nls}$ in Algorithm \ref{algorithm2} for solving \eqref{min_problem}. Define the perturbation energy function
\begin{equation*}
H\left( \bar{x},x \right) =f\left( \bar{x} \right) +g_1\left( \bar{x} \right) -g_2\left( \bar{x} \right) +\frac{L}{2}\left\| \bar{x}-x \right\| _{M}^{2}.
\end{equation*}
Suppose that $\left\{ x^n \right\} _{n}$ has a cluster $x^*$ and $H\left(\bar{x},x \right)$ is a KL function, then the following statements hold 
\begin{itemize}
\item[{\rm (i)}] \textbf{Monotonicity}: The sequence $\left\{ H\left(\bar{x}^{n},x^n \right) \right\} _{n > N_0}$ is monotonically decreasing and converges to some limit: $\lim_{n\rightarrow \infty} H\left( \bar{x}^{n},x^n \right) =\bar{\zeta}$.
\item[{\rm (ii)}] \textbf{Optimality}: Any cluster point of $\left\{ x^n \right\}_n $ is the critical point of the problem \eqref{min_problem}, i.e. it satisfies the first order optimality condition. Moreover, the iteration sequence $\{x^n\}_n$ is bounded.
\item[{\rm (iii)}] \textbf{Convergence}:  The sequence $\left\{ x^n \right\}_n $ is globally convergent with $\sum_{n=1}^{\infty}{\left\| x^{n}-x^{n-1} \right\|}<\infty $.
\end{itemize}
\end{theorem}

\begin{remark}
As for the framework of Algorithm \ref{algorithm2}, we can naturally obtain global convergence. Noting that the difference between Algorithm 1 and Algorithm \ref{algorithm2} is the difference of \eqref{convergence1} and \eqref{convergence2}. Therefore, in the framework of Algorithm \ref{algorithm2}, the subsequent proofs are not quite different from the proofs for Algorithm \ref{algorithm1}.
\end{remark}

\subsection{Convergence rate of the Algorithm \ref{algorithm1}}

Our analysis requires explicit characterization of the Kurdyka-\L ojasiewicz (KL) exponent associated with the auxiliary function $A(x,\bar{y},y)$ to quantify local convergence rates. The subsequent theorem, which adapts well-established convergence rate theory (cf. \cite[Theorem 2]{Attouch2009}, \cite[Lemma 1]{Artacho2018}), formalizes this relationship. 
\begin{theorem}[Local convergence rate]\label{thm:local:rate} 
Under the assumptions of Theorem \ref{thm:globalconvergence}, let $\{x^n\}_n$ be a sequence generated by Algorithm \ref{algorithm1} with a convergent subsequence $x^{n_i} \to x^*$. If $A(x,\bar{y},y)$ is a KL function with $\phi(s) = cs^{1-\theta}$ in the KL inequality $(\theta\in[0,1), c>0)$, then

\begin{itemize}
\item[{\rm (i)}] When $\theta = 0$, $\exists \  m_0 > 0$ such that 
      \begin{equation*}
        x^n = x^* \quad \forall \  n > m_0.
      \end{equation*}
\item[{\rm (ii)}] For $\theta \in (0, \frac{1}{2}]$, $\exists \  p_1 > 0$, $\eta \in (0,1)$, and $m_1 > 0$ such that 
      \begin{equation*}
        \|x^n - x^*\| \leq p_1\eta^{n} \quad \forall \  n > m_1.
      \end{equation*}
\item[{\rm (iii)}] For $\theta \in (\frac{1}{2}, 1)$, $\exists \ p_2 > 0$ and $m_2 > 0$ such that 
      \begin{equation*}
        \|x^n - x^*\| \leq p_2n^{-\frac{1-\theta}{2\theta-1}} \quad \forall \ n > m_2.
      \end{equation*}
\end{itemize}
\end{theorem}

\section{Numerical experiments}\label{sec:num}
In this section, we conduct two representative numerical experiments to validate our algorithm's advantages in solving non-convex optimization problems
\begin{itemize}
    \item \textbf{SCAD-Regularized Least Squares}: {Implemented an NVIDIA RTX 4070 GPU system (24GB VRAM).} 
    \item \textbf{Nonlocal Ginzburg-Landau Segmentation}: Executed on an NVIDIA RTX 3080 Ti GPU system (16GB VRAM).
\end{itemize}
Detailed implementation protocols and comparative results will be systematically presented in subsequent subsections. For all tables, when the iteration number exceeds the maximum iterations, we use the notation ``---" instead of recording the CPU time.

\begin{remark}\label{rem:KL:conver:two:models}
Before starting this subsection, we show that the minimization functions of the following two problems do have the KL property. The KL property of the original SCAD regularization \eqref{eq:lsp_l1} follows directly from \cite[Example 4.3]{Wen2018}. For the modified SCAD regularization in \eqref{eq:lsp_Huber}, since each $s_{M,i}(x_i)$ is a piecewise quadratic function, the KL property of $E(x)$ can be derived via the methodology in \cite[Section 5.2]{Li2018}.
Combining theoretical frameworks from \cite{Li2018,liu2019refined,ABS}, we confirm that both $A(x,\bar{y},y)$ and $H(x,w,\bar{y},y)$ are KL functions, thereby establishing global convergence. Finally, the discrete graphic Ginzburg-Landau functional, being a polynomial (hence semialgebraic) function, automatically satisfies the KL property \cite[Section 2.2]{ABS}.

\end{remark}

\subsection{Least squares problems with modified SCAD regularizer}\label{sec:num1}
We analyze a smoothed variant of the SCAD (Smoothly Clipped Absolute Deviation) regularizer through DC decomposition. Following the framework in \cite[Section 6.1]{Apx} and \cite{Wen2018}, the SCAD penalty admits
\begin{equation*}
S(x) = \mu\|x\|_1-\mu\sum_{i=1}^{k}\int_0^{|x_i|}\frac{[\min\{\theta\mu,u\}-\mu]_+}{(\theta-1)\mu}du,
\end{equation*}
where $\theta > 1$ controls concavity and $\mu>0$ regulates sparsity. The second term $\tilde S(x) =  \mu \|x\|_1 -S(x) = \sum_{i}^k \tilde s_i(x_i)$ decomposes into component functions 
$\tilde s_i(x_i)$ with piecewise structure
\begin{equation*}
    \tilde s_i(x_i) = \left\{
        \begin{array}{ll}
         0& |x_i|\leq\mu \\
         \frac{(|x_i| - \mu)^2}{2(\theta - 1)\mu}& \mu<|x_i|<\theta \\
         \mu|x_i| - \frac{\mu^2(\theta+1)}{2}& |x_i|\geq\theta\mu 
    \end{array},
        \right. \ \ \nabla_i \tilde S_i(x_i)= \text{sign}(x_i) \dfrac{[\min\{\theta \mu,|x_i|\} -\mu]_{+}} {(\theta-1)}.
\end{equation*}
We formulate the $l_1$ norm regularized least squares problem with the SCAD regularizer as follows
\begin{equation}\label{eq:lsp_l1}
    \min_{x\in \mathbb{R}^k}E(x) = \frac12\|Ax-b\|^2 + S(x).
\end{equation}
To ensure differentiability, we replace the $l_1$ norm with the Huber function $\mathcal{H}_{\alpha}(x):= \sum_{i}\mathcal{H}(x_i,\alpha)$
\begin{equation*}
    \mathcal{H}(x_i,\alpha) = \left\{
    \begin{array}{ll}
         \frac{|x_i|^2}{2\alpha}& |x_i|\leq\alpha \\
         |x_i|-\frac{\alpha}{2}&  |x_i|>\alpha
    \end{array},
    \right. 
\end{equation*}
where $\alpha=0.5\mu$ balances robustness. The modified SCAD-regularized least squares problem becomes
\begin{equation}\label{eq:lsp_Huber}
    \min_{x\in \mathbb{R}^k}E(x) = \frac12\|Ax-b\|^2 + \mu \mathcal{H}_{\alpha}(x) - \tilde S(x),
\end{equation}
with $A\in \mathbb{R}^{m\times k}$, $b\in \mathbb{R}^{m}$. The composite penalty  $S_{M}(x): =  \mu\mathcal{H}_{\alpha}(x) - \tilde S(x) := \sum_{i=1}^k s_{M,i}(x_i)$, where
   $p_{M,i}(x_i)$ exhibits piecewise behavior
        \begin{equation}\label{eq:huber-scad}
        \frac{s_{M,i}(x_i)}{\mu} = \left\{
        \begin{array}{ll}
       |x_i|^2/(2\alpha)& |x_i|\leq\alpha \\
         |x_i|-{\alpha}/{2}&  \alpha < |x_i|\leq\mu \\
         |x_i|-{\alpha}/{2} -(|x_i| - \mu)^2/(2(\theta - 1)\mu)& \mu<|x_i|<\theta\mu \\
         (\mu(\theta+1) - \alpha)/{2} & |x_i|\geq\theta\mu 
    \end{array}.
        \right.
    \end{equation}
    
We validate on LIBSVM datasets \cite{chih2011libsvm} converted to MATLAB format for binary classification. Initialized at $x^0=0$, convergence is triggered when
\begin{align*}
\frac{\|x^n - x^{n-1}\|}{\max\{1,\|x^n\|\}}< 
\epsilon \text{, where } \epsilon = 10^{-i} \text{, for }i = 4,5, \ldots, 9.    
\end{align*}

We will compute the original SCAD \eqref{eq:lsp_l1} and the modified SCAD \eqref{eq:lsp_Huber} for comparisons. From \cite{Shensun2023}, we see that the Huber-SCAD model \eqref{eq:lsp_Huber} can obtain nearly the same sparsity as the original SCAD model. We refer to \cite{Li_2020} for other types of DC-type regularization for signal and image processing. 

Evaluating various algorithms involves recording both iteration counts (iter) and CPU time. We present computational results in Table \ref{table:l1} and Table \ref{table:Huber},  corresponding to problem \eqref{eq:lsp_l1} and problem \eqref{eq:lsp_Huber} with parameters $\mu = 5\times10^{-4}$ and $\theta = 10$, respectively. In both numerical experiments, we used LSDE \eqref{LSDE} to update the extrapolation parameter $\beta_n$. We compare five algorithms to solve these two problems: our proposed algorithm npDCAe$_{\text{nls}}$ (pDCAe$_{\text{nls}}$), pDCAe, pDCAe (no-restart), along with the DC algorithm (DCA) and the BDCA-Backtracking algorithm (BDCA$_{\text{ls}}$) from \cite{Artacho2018}. Due to \eqref{rem:6} in the following, we discuss the details of the implementation of these algorithms only with Algorithm \ref{algorithm1}.

\begin{remark}\label{rem:6}
It is worth noting that in solving \eqref{eq:lsp_l1} and \eqref{eq:lsp_Huber}, Algorithm \ref{algorithm1} and Algorithm \ref{algorithm2} become equivalent when appropriate preconditioners are selected. Take \eqref{eq:lsp_l1} as an example, when we set $f(x)$, $g_1(x)$ and $g_2(x)$ as follows:
\[
f(x)=\frac12\|Ax-b\|^2, \  g_1(x)=\mu\|x\|_1 \ \text{ and  } \ g_2(x)=\mu\sum_{i=1}^{k}\int_0^{|x_i|}\frac{[\min\{\theta\mu,u\}-\mu]_+}{(\theta-1)\mu}du. 
\]
Choosing $M=\lambda_{A^TA}I-A^TA$, we can see that the two algorithms have the same iterations.
\end{remark}

\begin{itemize}
    \item \textbf{pDCAe}:
    For this algorithm, we set $L=\lambda_{A^TA}$, where $\lambda_{A^TA}$ represents the largest eigenvalue of the matrix $A^TA$, and choose the extrapolation parameter $\left\{ \beta _n \right\}$ as FISTA, and perform both fixed restart (with $\bar{T}=200$) and adaptive restart strategy as
 described above \eqref{convention_restart}. One can refer to   \cite[Section 5.1]{Wen2018} for more details.
\item \textbf{pDCAe(no-restart)}: This special version of the algorithm {pDCAe} does not perform the fixed restart or the adaptive restart strategy. The extrapolation parameter of this algorithm has been chosen as FISTA. This comparison aims to show the importance of the restart of the extrapolation parameter.
\item \textbf{npDCAe$_{\text{nls}}$}(\textbf{pDCAe$_{\text{nls}}$}): This algorithm combined pDCAe   \cite{Wen2018} and a non-monotone line search technique, together with LSDE \eqref{LSDE} as detailed in Algorithm \ref{algorithm1}. The parameters used by this algorithm to solve \eqref{eq:lsp_Huber} are chosen as follows: $\lambda_{\text{max}} = 2$, $\alpha = 10^{-8}$, $N_{\text{max}} = 3$, $\rho = 0.3$, $\omega = 0.9$, $\eta 
 = 1.9$ ($\omega$ and $\eta$ may vary depending on accuracy), $b_1 = 0.001$, $b_2=0$.
\item \textbf{BDCA$_{\text{ls}}$}: This algorithm is based on a combination of DCA together with a line search technique, as detailed in    \cite[Algorithm 2]{Artacho2018}. To simplify computations, we employ distinct convex splitting to sidestep equation solving. The convex splitting can be defined as follows
    \begin{equation} \label{eq:another:split}
        E(x) =\mu \mathcal{H}_{\alpha}(x) + \frac{\lambda_{A^TA}}{2}\|x\|^2 - \left(\frac{\lambda_{A^TA}}{2}\|x\|^2 + \tilde S(x) - \frac12\|Ax-b\|^2\right),
    \end{equation}
    where $\lambda_{A^TA}$ represents the maximum eigenvalue of the matrix $A^TA$.
    \item \textbf{DCA}: This is the classical DC algorithm, which is a special version of the algorithm BDCA$_{\text{ls}}$ without the line search    \cite[Algorithm 1]{Artacho2018}. 
\end{itemize}

\begin{remark}\label{rem:cri:new}
In addition, we introduce an additional criterion to demonstrate the accuracy of the first-order optimality condition of the problem \eqref{eq:lsp_l1} and \eqref{eq:lsp_Huber}. For  \eqref{eq:lsp_l1}, the following criterion is used to measure the residual of the first-order optimality condition for the nonsmooth problem \eqref{eq:lsp_l1} 
\begin{equation*}
res\left( x^n \right) =\left\| G_{\proj}\left( x^n \right) \right\| ,\quad G_{\proj}\left( x^n \right) =\left( x^n-\prox_{I_2}\left( x^n-\nabla I_1\left( x^n \right) \right) \right),
\end{equation*}
where $I_1\left( x \right) =\frac{1}{2}\left\| Ax-b \right\| ^2+p_2\left( x \right) $,  
$I_2\left( x \right) =\mu \left\| x \right\| _1$,  and  $\prox_{I_2}(z):= \argmin_{x} \{I_2(x)+\frac{1}{2}\|x-z\|^2$\}.
For \eqref{eq:lsp_Huber}, we choose $res\left( x^n \right):=\|\nabla E(x^n)\|$ with $E(x)$ as in \eqref{eq:lsp_Huber}.
\end{remark}

\begin{table}[htbp!]
\centering
\captionsetup{justification = centering}
\caption{Solving (\ref{eq:lsp_l1}) on random instances\\ (\textcircled{1}pDCAe 
\textcircled{2}pDCAe(no-restart)
\textcircled{3}DCA  
\textcircled{4}BDCA$_{\text{ls}}$
\textcircled{5}npDCAe$_{\text{nls}}$, maxier(Max)=5e6)}
\scalebox{0.9}{
\label{table:l1} 
\begin{tabular}{c|rrrrr|rrrrr}
\hline
\multicolumn{1}{c|}{Precision} & \multicolumn{5}{c|}{Iter} & \multicolumn{5}{c}{GPU time (s)} \\ \hline
Algorithm &\textcircled{1}  &
\textcircled{2} &
\textcircled{3} & \textcircled{4}  & \textcircled{5}  &\textcircled{1}  &
\textcircled{2} &
\textcircled{3} & \textcircled{4}   & \textcircled{5} \\ \hline
{$10^{-4}$} &  202 & 35851 & 1237  &  189  &  \textbf{74}   &   0.3  & 29.3 &   1.3  &   0.7  &  \textbf{0.2} \\
{$10^{-5}$} &  1002 &  139648  &  20644 & 5310 &   \textbf{586}   &        1.5  & 114.7 & 16.4  &  18.5 &   \textbf{1.5} \\
{$10^{-6}$} &  7802 &  Max & 169124 &  23630 &  \textbf{2610}   &           7.7  &  --- & 133.5  &  62.9 &  \textbf{6.4}  \\
{$10^{-7}$} &  39202 & Max & Max  &  138017 &  \textbf{7827}   &           37.6  &  --- & --- &  301.1 &  \textbf{18.1} \\
{$10^{-8}$} &  130202 &  Max & Max  &  Max &  \textbf{10866}   &           128.7  & --- & --- &  --- &  \textbf{23.1} \\
{$10^{-9}$} &  440402 &  Max & Max & Max &   \textbf{14072}   &          399.4  & --- & ---  & --- &  \textbf{31.2} \\ \hline
\end{tabular}}

\end{table}

\begin{table}[htbp!]
\centering
\captionsetup{justification = centering}
\caption{Solving (\ref{eq:lsp_Huber}) on random instances\\ (\textcircled{1}pDCAe 
\textcircled{2}pDCAe(no-restart)
\textcircled{3}DCA  
\textcircled{4}BDCA$_{\text{ls}}$
\textcircled{5}npDCAe$_{\text{nls}}$, maxier(Max)=5e6)}
\scalebox{0.9}{
\label{table:Huber}
\begin{tabular}{c|rrrrr|rrrrrrr}
\hline
\multicolumn{1}{c|}{Precision} & \multicolumn{5}{c|}{Iter} & \multicolumn{5}{c}{GPU time (s)} \\ \hline
Algorithm &\textcircled{1}  &
\textcircled{2} &
\textcircled{3} & 
\textcircled{4}  & 
\textcircled{5}  &
\textcircled{1}  &
\textcircled{2} &
\textcircled{3} & 
\textcircled{4}   & 
\textcircled{5} \\ \hline
{$10^{-4}$} &  202 & 42804 & 1216  &  187  &  \textbf{100}   &   0.5  & 35.0 &   1.4  &   0.6  &  \textbf{0.3} \\
{$10^{-5}$} &  1002 &  161030  &  20132 & 4243 &   \textbf{506}   &        1.5  & 131.9 & 16.8  &  16.1 &  \textbf{1.5} \\
{$10^{-6}$} &  8002 &  Max & 163751 &  28232 &  \textbf{2805}   &           8.2  &  --- & 144.2  &  98.2 &  \textbf{7.8}  \\
{$10^{-7}$} &  38402 & Max & Max  &  195804 &  \textbf{9672}   &           37.6  &  --- & --- &  593.4 &  \textbf{26.7} \\
{$10^{-8}$} &  132802 &  Max & Max  &  Max &  \textbf{13545}   &           127.8  & --- & --- &  --- &  \textbf{35.4} \\
{$10^{-9}$} &  437402 &  Max & Max & Max &   \textbf{17067}   &          473.2  & --- & ---  & --- &  \textbf{50.2} \\ \hline
\end{tabular}}
\end{table}

\begin{figure}[htbp!]   
  \centering            
   \subfloat[$l_1$-Convergence-Rate (GPU)]{\label{fig:l1_convergence_rate_GPU}\includegraphics[width=0.45\textwidth]{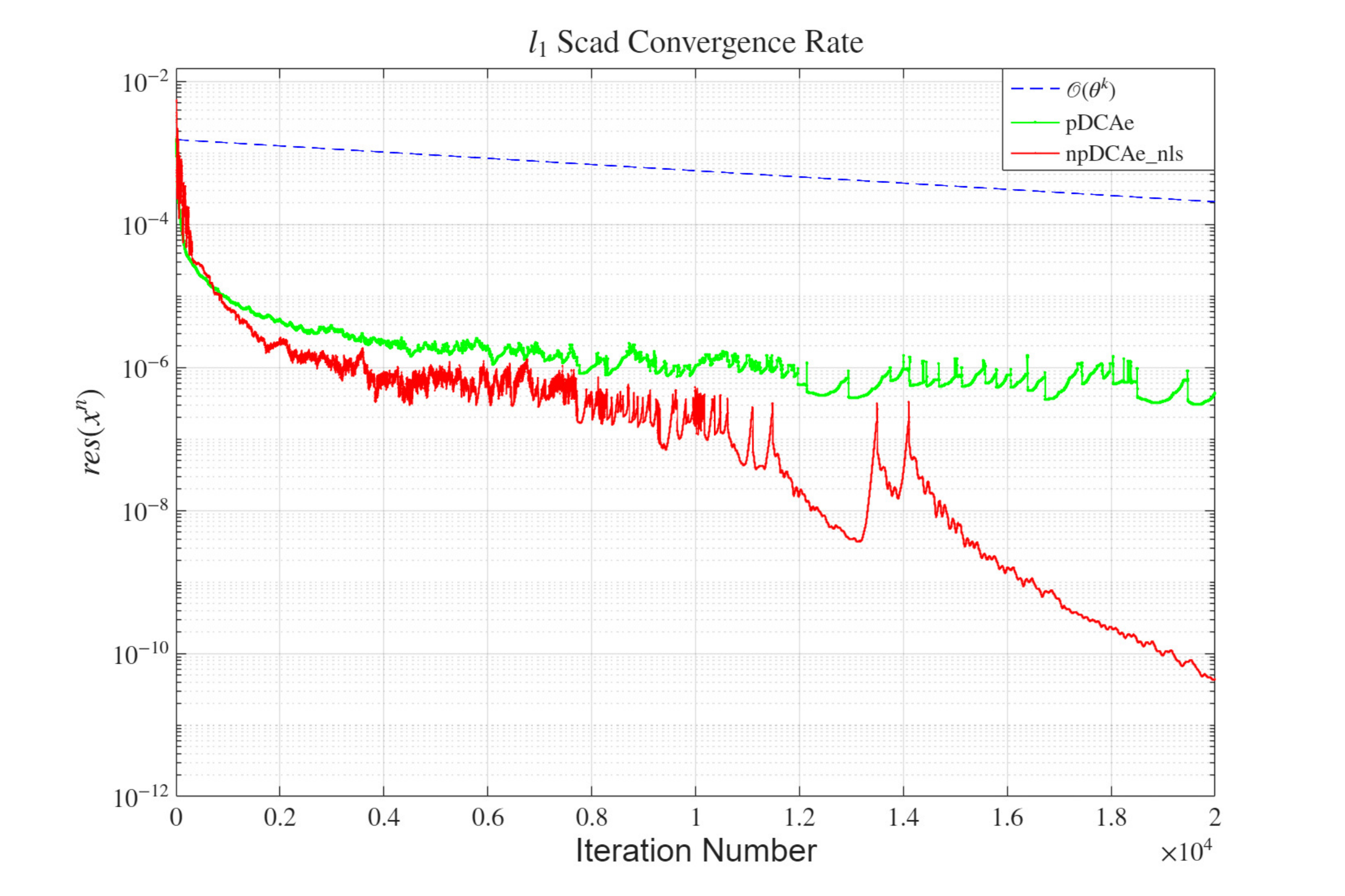}}\quad \quad 
   \subfloat[Huber-Convergence-Rate (GPU)]{\label{fig:Huber_convergence_rate_GPU}\includegraphics[width=0.45\textwidth]{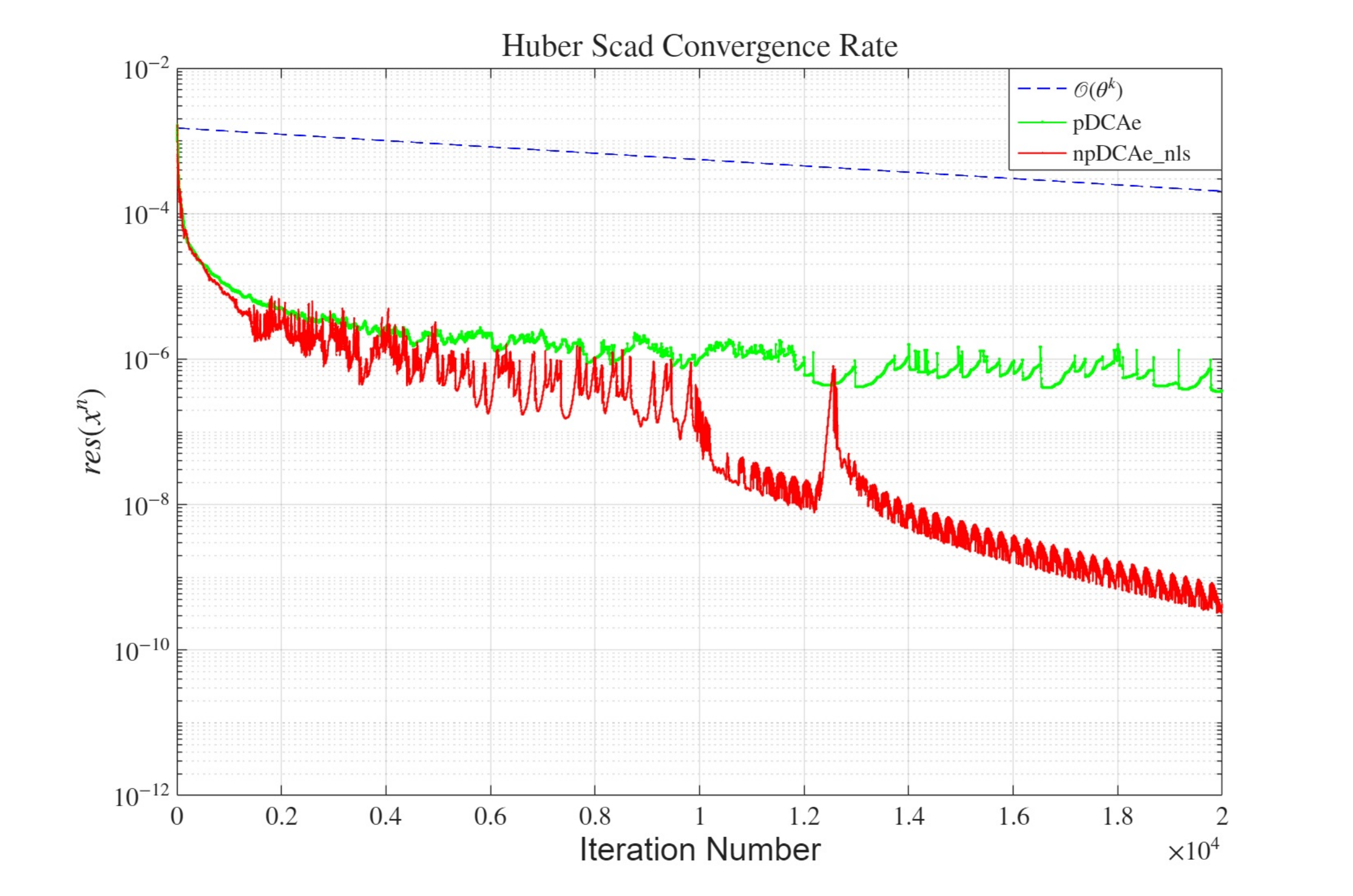}}\quad
  \caption{The convergence rate of two algorithms for solving \eqref{eq:lsp_l1} and \eqref{eq:lsp_Huber}.}\label{fig:rate:scad}       
\end{figure}
 Table \ref{table:l1} and Table \ref{table:Huber}  are the numerical comparisons for solving the SCAD model \eqref{eq:lsp_l1} and the modified SCAD model \eqref{eq:lsp_Huber}, respectively. It can be seen that the proposed npDCAe$_{\text{nls}}$ is highly efficient compared to pDCAe. Especially, npDCAe$_{\text{nls}}$ is much faster for very high accuracy, i.e., $\|x^{n+1}-x^n\|<10^{-8}$. Figure \ref{fig:rate:scad} shows that the proposed algorithm can obtain a promising linear convergence rate numerically for both the SCAD model and the modified SCAD model.
From \cite[Theorem 4.2]{liu2019refined} we know that the KL exponent of the minimization functions \eqref{eq:lsp_l1} and \eqref{eq:lsp_Huber} are both $\frac{1}{2}$. In the Figures \eqref{fig:l1_convergence_rate_GPU} and \eqref{fig:Huber_convergence_rate_GPU}, we can see the fast local convergence, which can be faster than linear convergence. 


\subsection{Graphic Ginzburg-Landau model for image segmentation} \label{sec:num2}
In this subsection, we focus on a segmentation problem using a nonlocal Ginzburg-Landau model with prior guidance, formulated as \cite{BF}
\begin{equation}\label{eq:nonlocalGL}
    \min_{x\in \mathbb{R}^N}E(x) = \sum_{i,j}\frac{\tau}{2}w_{ij}(x(i)-x(j))^2 +\frac{1}{\tau}\bb{W}(x) + \frac{\gamma}{2}\sum_i\Lambda(i)(x(i)-y(i))^2,
\end{equation}
where $x$ denotes segmentation labels, $w_{ij} =K(i,j) \cdot N(i,j)$ encodes nonlocal interactions via Gaussian similarity $K(i,j) = \exp\left(-\|P_{i} - P_{j}\|_2^2/\kappa^2\right)$ and binary adjacency $N(i,j)$, $\mathbb{W}(x) = \frac14\sum_{i=1}^{N}(x(i)^2 - 1)^2$ enforces binary constraints, the diagonal matrix $\Lambda$ and the vector $y$ incorporate prior labels \cite{Shensun2023}. We choose a $25 \times 25$ square box surrounding each pixel for nonlocal interactions. Parameters $\tau = \gamma = 10$ balance regularization strength. 

The energy function mentioned above $E(x)$ is divided into two components
\begin{align*}
&f(x) = \frac{\gamma}{2}\sum_i\Lambda(i)(x(i)-y(i))^2, \ g_1(x) = \sum_{ij}\frac{\tau}{2}w_{ij}(x(i)-x(j))^2  \text{ and } g_2(x) = \frac{1}{\tau}\mathbb{W}(x).     
\end{align*}
For comparison, we implement the following algorithms
\begin{itemize}
    \item \textbf{Preconditioned DC algorithms}: pDCA \cite{LeThi2018}, pDCAe \cite[Section 3]{Wen2018}, pDCA$_{\text{nls}}$ \cite[Section 3]{Artacho2018}, pDCAe$_{\text{nls}}$, npDCAe \cite{DS},
    npDCAe$_{\text{nls}}$.
    \item \textbf{CG(conjugate gradient)-based methods}: pDCAe(CG) and pDCAe$_{\text{nls}}$(CG).
\end{itemize}
For all the preconditioned DC algorithms, we use 5 Jacobi iterations, while for pDCAe(CG) and pDCAe$_{\text{nls}}$(CG), we employ CG with tolerance $\|\bar{x}^{l+1} - \bar{x}^l\|<10^{-11}$. For pDCAe$_{\text{nls}}$, npDCAe$_{\text{nls}}$ and pDCAe$_{\text{nls}}$(CG), we employ the following parameters: $\lambda_{\max}=2$, $N_{\max}=3$, $\rho = 0.3$, $\omega=0.001$, $\eta =0.3$, $b_1=0.001$, $b_2=0$.

The first criterion ``DICE Bound" in Table \ref{tab:nonlocalGLmodel}, relies on the segmentation coefficient DICE, a metric used to evaluate the segmentation outcome. It is calculated using the formula $\text{DICE} = \frac{2|X\cap Y|}{|X|+|Y|}$, where $|X|$ and $|Y|$ denote the pixel counts in the segmentation and ground truth, respectively. Moreover, the criterion ``DICE Bound" refers to the DICE coefficient of the segmentation, achieving a value of $0.985$, indicating a high degree of similarity.

\begin{figure}[htbp]   
  \centering            
  \subfloat[Clover]   
  {\label{fig:clover1}\includegraphics[width=0.18\textwidth]{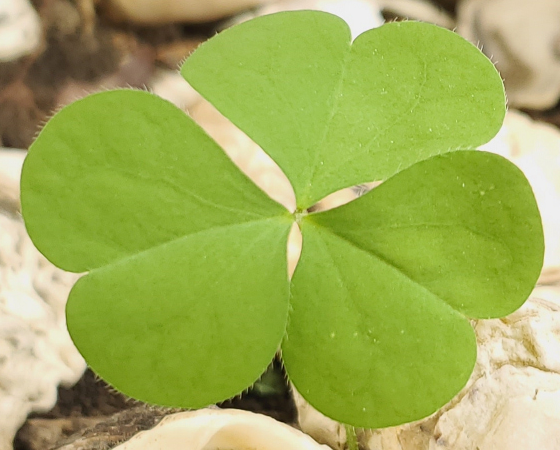}}\quad
  \subfloat[Label]
  {\label{fig:cloverlabel}\includegraphics[width=0.18\textwidth]{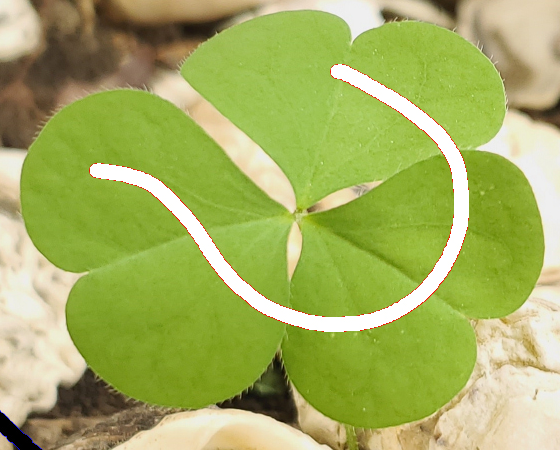}}\quad
  \subfloat[result 1]
  {\label{fig:result1}\includegraphics[width=0.18\textwidth]{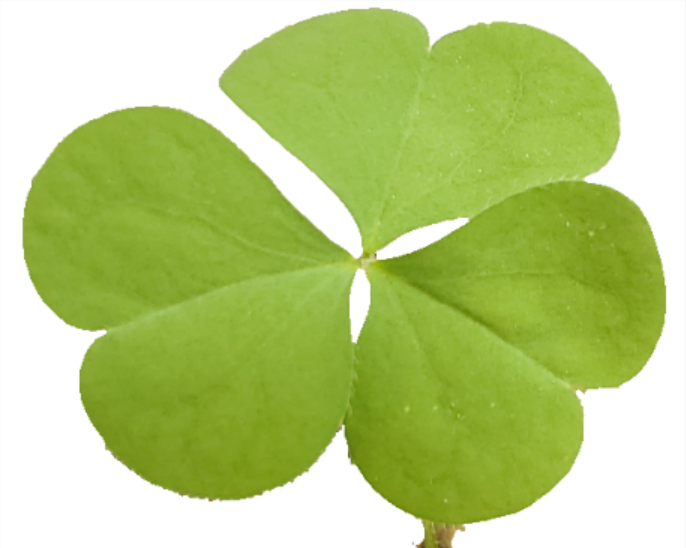}}\quad
  \subfloat[result 2]
  {\label{fig:result2}\includegraphics[width=0.18\textwidth]{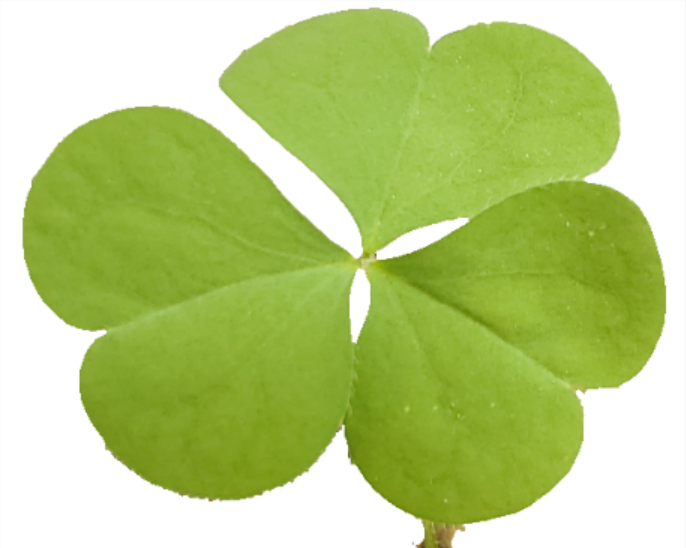}}\quad
  \subfloat[result 3]
  {\label{fig:result3}\includegraphics[width=0.18\textwidth]{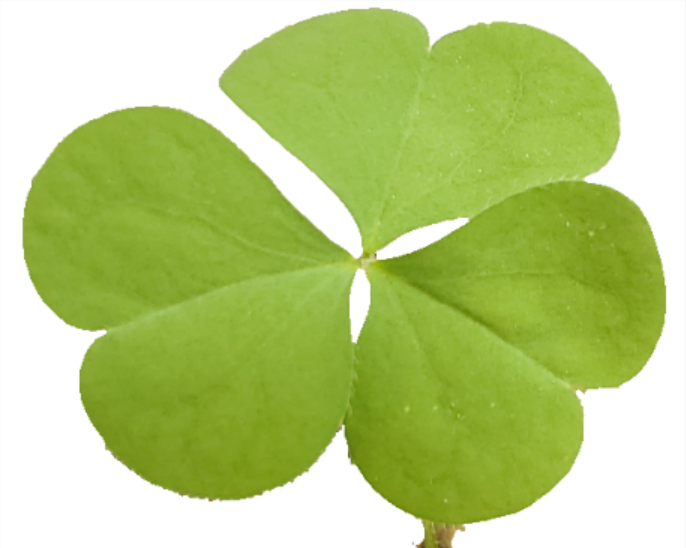}}\quad
  \caption{The performance of segmentation assignment.}\label{fig:clover}
  \label{fig:clover_seg}         
\end{figure}

In Figure \ref{fig:clover_seg}, we present the segmentation results for an illustrative case: figure \ref{fig:clover1} is the original input image, figure \ref{fig:cloverlabel} labels the prior that the white pixels label the positive prior, and the black
 pixels label the negative prior, while the remaining figures show the segmentation results under
 three different criteria
(DICE Bound, $\|\nabla E(x)\|\leq 10^{-11}$ and $\|x^n-x^{n-1}\|\leq 10^{-11}$).

From Table \ref{tab:nonlocalGLmodel}, it can be seen that both pDCAe$_{\text{nls}}$ and npDCAe$_{\text{nls}}$ can be 30\%-40\% faster than their pDCAe counterparts. For the image segmentation task, it is also highly efficient to obtain a high DICE with a GPU.

From Tables \ref{table:l1}, \ref{table:Huber} and \ref{tab:nonlocalGLmodel}, we use LSDE \eqref{LSDE} to update the extrapolation parameter $\beta_n$, it can be seen that the algorithms pDCAe$_{\text{nls}}$ and npDCAe$_{\text{nls}}$ outperform other algorithms in both iterations numbers, requiring the least amount of time to fulfill the termination criteria.


\begin{table}[htbp!]
\centering
\captionsetup{justification = centering}
\caption{Solving \ref{eq:nonlocalGL} on two different termination criteria.(Criteria \uppercase\expandafter{\romannumeral1}: $\|x^n-x^{n-1}\|$, Criteria \uppercase\expandafter{\romannumeral2}: $\|\nabla E(x)\|$, maxiter(Max)=5000)
\textcircled{1}pDCAe(CG)
\textcircled{2}pDCAe$_{\text{nls}}$\text{(CG)}
\textcircled{3}pDCA
\textcircled{4}pDCA$_{\text{nls}}$
\textcircled{5}pDCAe 
\textcircled{6}pDCAe$_{\text{nls}}$
\textcircled{7}npDCAe
\textcircled{8}npDCAe$_{\text{nls}}$
}
\label{tab:nonlocalGLmodel}
 \scalebox{0.9}{
\begin{tabular}{cccrrrrrrrrr}
\hline
\multicolumn{3}{c}{Algorithm} &{\textcircled{1}} &
{\textcircled{2}} & 
{\textcircled{3}} &
{\textcircled{4}} & 
{\textcircled{5}} &
{\textcircled{6}} &
{\textcircled{7}} &
{\textcircled{8}} & \\ \hline
\multicolumn{2}{c}{\multirow{2}{*}{DICE Bound}} & Iter & 88 & \textbf{54} & 1003 & 341 & 191 & \textbf{53} & 58 & \textbf{39} \\ 
& &Time(s)& 45.67 & \textbf{29.34} & 76.72 & 38.12& 19.81 & \textbf{7.72} & 8.56 & \textbf{7.83} \\ \hline
\multirow{6}{*}{\uppercase\expandafter{\romannumeral1}}& \multirow{2}{*}{$10^{-1}$}& Iter & 119 & \textbf{80} & 1052 & 442 & 238 & \textbf{72} & 160 & \textbf{55} \\
& &Time(s)& 59.67 & \textbf{41.12}& 80.72 & 43.12& 21.81 & \textbf{10.72} & 16.06 & \textbf{9.43} \\
 & \multirow{2}{*}{$10^{-3}$}&Iter & 315 & \textbf{174} & 3182 & 1284 & 739 & \textbf{231} & 291 & \textbf{115} \\
& &Time(s)& 137.48 & \textbf{91.36} & 237.24 & 126.04 & 65.40 & \textbf{23.86} & 26.33 & \textbf{13.54} \\ 
& \multirow{2}{*}{$10^{-5}$} & Iter & 520 & \textbf{290} & 4891 & 2154 & 1201 & \textbf{399} & 411 &\textbf{197} \\
& &Time(s)& 218.65 & \textbf{151.98} & 92.3 & 192.1 & 60.6 & \textbf{41.18} & 65.5 & \textbf{21.19} \\  
& \multirow{2}{*}{$10^{-7}$} & Iter & 720 & \textbf{462} & Max & 3277 & 1801 & \textbf{550} & 554 &\textbf{272} \\
& &Time(s)& 298.64 & \textbf{241.24} & --- & 291.58 & 153.98 & \textbf{60.65} & 50.06 & \textbf{27.84} \\  
& \multirow{2}{*}{$10^{-9}$} & Iter & 1001 & \textbf{600} & Max & 4641 & 2601 & \textbf{987} & 741 &\textbf{395} \\
& &Time(s)& 416.17 & \textbf{316.12} & --- & 408.92 & 216.26 & \textbf{106.15} & 68.41 & \textbf{38.45} \\  
& \multirow{2}{*}{$10^{-11}$} & Iter & 1320 & \textbf{655} & Max & Max & 3201 & \textbf{1203} & 931 &\textbf{589} \\
& &Time(s)& 537.25 & \textbf{345.42} & --- & --- & 254.27 & \textbf{126.37} & 86.24 & \textbf{54.79} \\  \hline
\multirow{6}{*}{\uppercase\expandafter{\romannumeral2}} & \multirow{2}{*}{$10^{-1}$} &Iter& 136 & \textbf{78} & 1188 & 485 & 270 & \textbf{82} & 165 & \textbf{85} \\
& &Time(s)& 67.52 & \textbf{40.73} & 89.09 & 46.86 & 24.36 & \textbf{11.28} & 16.47 & \textbf{11.66} \\ 
& \multirow{2}{*}{$10^{-3}$}  &Iter& 325 & \textbf{167} & 3535 & 1393 & 771 & \textbf{253} & 285 & \textbf{114} \\
& &Time(s)& 141.35 & \textbf{87.76} & 237.24 & 126.04 & 65.4 & \textbf{26.86} & 26.33 & \textbf{13.89} \\ 
& \multirow{2}{*}{$10^{-5}$}  &Iter& 562 & \textbf{271} & Max & 2672 & 1399 & \textbf{409} & 405 & \textbf{198} \\
& &Time(s)& 235.65 & \textbf{141.43} & --- & 239.17 & 118.97 & \textbf{42.78} & 36.33 & \textbf{22.21} \\ 
& \multirow{2}{*}{$10^{-7}$}  &Iter& 849 & \textbf{409} & Max & 3445 & 2117 & \textbf{755} & 497 & \textbf{278} \\
& &Time(s)& 350.93 & \textbf{213.45} & --- & 306.12 & 181.71 & \textbf{79.58} & 46.47 & \textbf{28.31} \\ 
& \multirow{2}{*}{$10^{-9}$}  &Iter& 1099 & \textbf{556} & Max & 4656 & 2819 & \textbf{1084} & 687 & \textbf{356} \\
& &Time(s)& 449.59 & \textbf{291.34} & --- & 410.21 & 230.09 & \textbf{115.96} & 62.01 & \textbf{35.39} \\ 
& \multirow{2}{*}{$10^{-11}$}&Iter&  1379 & \textbf{670} & Max & Max & 3525 & \textbf{1232} & 887 & \textbf{508} \\
& &Time(s)& 561.77 & \textbf{353.76} & --- & --- & 274.84 & \textbf{131.98} & 81.72 & \textbf{48.82} \\ \hline
\end{tabular}
}
\end{table}

\section{Conclusion}\label{sec:conclu}
We investigate two kinds of proximal DC algorithms and their corresponding extrapolation, which are commonly used in non-convex optimization. The extrapolation parameter is completely determined by a non-monotone line search method. We prove the global convergence of the improved algorithm with the KL property. The numerical experiments show that the proposed algorithms are highly efficient. In addition, our non-monotone line search can also be performed using both methods from \cite{Dai2002} and \cite{lu2019nonmonotone}, which are also discussed in   \cite{ferreira2024boosted}.

\noindent
{\small
	\textbf{Acknowledgements}
Ran Zhang and Hongpeng Sun acknowledge the support of National Key R\&D Program of China (2022ZD0116800), the National Natural Science Foundation of China under grant No. \,12271521, and Beijing Natural Science Foundation No. Z210001. The second author acknowledges Prof. Ting Kei Pong of The Hong Kong Polytechnic University for suggesting the techniques developed in \cite{liu2019refined} during a private communication. 
}

\bibliographystyle{plain}
\bibliography{bdfab_ls}

\begin{thebibliography}{10}

\bibitem{Apx}
Miju Ahn, Jong-Shi Pang, and Jack Xin.
\newblock Difference-of-convex learning: Directional stationarity, optimality, and sparsity.
\newblock {\em SIAM J. Optim.}, 27(3):1637--1665, 2017.

\bibitem{Artacho2018}
Francisco~J. Arag{\'o}n~Artacho, Ronan M.~T. Fleming, and Phan~T. Vuong.
\newblock Accelerating the dc algorithm for smooth functions.
\newblock {\em Math. Program.}, 169(1):95--118, May 2018.

\bibitem{ArtachoSIAM}
Francisco~J. Arag\'{o}n~Artacho and Phan~T. Vuong.
\newblock The boosted difference of convex functions algorithm for nonsmooth functions.
\newblock {\em SIAM J. Optim.}, 30(1):980--1006, 2020.

\bibitem{ABS}
H.~Attouch, J.~Bolte, and B.~F. Svaiter.
\newblock Convergence of descent methods for semi-algebraic and tame problems: proximal algorithms, forward--backward splitting, and regularized gauss--seidel methods.
\newblock {\em Math. Program.}, 137(1):91--129, Feb 2013.

\bibitem{Attouch2009}
Hedy Attouch and J{\'e}r{\^o}me Bolte.
\newblock On the convergence of the proximal algorithm for nonsmooth functions involving analytic features.
\newblock {\em Math. Program.}, 116(1):5--16, Jan 2009.

\bibitem{aujol2024parameter}
Jean-Fran{\c{c}}ois Aujol, Luca Calatroni, Charles Dossal, Hippolyte Labarri{\`e}re, and Aude Rondepierre.
\newblock Parameter-free fista by adaptive restart and backtracking.
\newblock {\em SIAM J. Optim.}, 34(4):3259--3285, 2024.

\bibitem{HBPL}
H.~H. Bauschke and P.~L. Combettes.
\newblock {\em Convex Analysis and Monotone Operator Theory in Hilbert Spaces}.
\newblock Springer Cham, seond edition, 2017.

\bibitem{BF}
A.~L. Bertozzi and A.~Flenner.
\newblock Diffuse interface models on graphs for classification of high dimensional data.
\newblock {\em SIAM Review}, 58(2):293--328, 2016.

\bibitem{BSCC}
Kristian Bredies and Hongpeng Sun.
\newblock Preconditioned douglas--rachford splitting methods for convex-concave saddle-point problems.
\newblock {\em SIAM J. Numer. Anal.}, 53(1):421--444, 2015.

\bibitem{chih2011libsvm}
Chang Chih-Chung.
\newblock Libsvm: a library for support vector machines.
\newblock {\em ACM transactions on intelligent systems and technology}, 2:27--1, 2011.

\bibitem{Pangcui}
Ying Cui and Jong-Shi Pang.
\newblock {\em Modern Nonconvex Nondifferentiable Optimization}.
\newblock Society for Industrial and Applied Mathematics, Philadelphia, PA, 2021.

\bibitem{Dai2002}
Y.~H. Dai.
\newblock On the nonmonotone line search.
\newblock {\em J. Optim. Theory Appl.}, 112(2):315--330, 2002.

\bibitem{DS}
S.~Deng and H.~Sun.
\newblock A preconditioned difference of convex algorithm for truncated quadratic regularization with application to imaging.
\newblock {\em J. Sci. Comput.}, 88(2):1--28, 2021.

\bibitem{fercoq2019adaptive}
Olivier Fercoq and Zheng Qu.
\newblock Adaptive restart of accelerated gradient methods under local quadratic growth condition.
\newblock {\em IMA J. Numer. Anal.}, 39(4):2069--2095, 2019.

\bibitem{ferreira2024boosted}
Orizon~P Ferreira, Elianderson~M Santos, and Jo{\~a}o Carlos~O Souza.
\newblock A boosted dc algorithm for non-differentiable dc components with non-monotone line search.
\newblock {\em Comput. Optim. Appl.}, pages 1--36, 2024.

\bibitem{hinder2020generic}
Oliver Hinder and Miles Lubin.
\newblock A generic adaptive restart scheme with applications to saddle point algorithms.
\newblock {\em arXiv preprint arXiv:2006.08484}, 2020.

\bibitem{LPH}
H.~A. Le~Thi and D.~T. Pham.
\newblock Convex analysis approach to d.c. programming: theory, algorithms and applications.
\newblock {\em Acta Mthematics Vietnamica}, 22(1):289--355, 1997.

\bibitem{LeThi2018}
Hoai~An Le~Thi and Tao Pham~Dinh.
\newblock Dc programming and dca: thirty years of developments.
\newblock {\em Math. Program.}, 169(1):5--68, May 2018.

\bibitem{Lethi2024}
Hoai~An Le~Thi and Tao Pham~Dinh.
\newblock Open issues and recent advances in dc programming and dca.
\newblock {\em J. Global Optim.}, 88(3):533--590, Mar 2024.

\bibitem{Li2018}
Guoyin Li and Ting~Kei Pong.
\newblock Calculus of the exponent of kurdyka--{\l}ojasiewicz inequality and its applications to linear convergence of first-order methods.
\newblock {\em Found. Comput. Math.}, 18(5):1199--1232, Oct 2018.

\bibitem{Li_2020}
Peng Li, Wengu Chen, Huanmin Ge, and Michael~K Ng.
\newblock $l_1$ $-$ $\alpha l_2$ minimization methods for signal and image reconstruction with impulsive noise removal.
\newblock {\em Inverse Problems}, 36(5):055009, apr 2020.

\bibitem{liu2019refined}
Tianxiang Liu, Ting~Kei Pong, and Akiko Takeda.
\newblock A refined convergence analysis of pdca e with applications to simultaneous sparse recovery and outlier detection.
\newblock {\em Comput. Optim. Appl.}, 73(1):69--100, 2019.

\bibitem{lu2019nonmonotone}
Zhaosong Lu and Zirui Zhou.
\newblock Nonmonotone enhanced proximal dc algorithms for a class of structured nonsmooth dc programming.
\newblock {\em SIAM J. Optim.}, 29(4):2725--2752, 2019.

\bibitem{o2015adaptive}
Brendan O’donoghue and Emmanuel Candes.
\newblock Adaptive restart for accelerated gradient schemes.
\newblock {\em Found. Comput. Math.}, 15:715--732, 2015.

\bibitem{ijcai2018p190}
Duy~Nhat Phan, Hoai~Minh Le, and Hoai An~Le Thi.
\newblock Accelerated difference of convex functions algorithm and its application to sparse binary logistic regression.
\newblock In {\em IJCAI-18}, pages 1369--1375, 2018.

\bibitem{Roc1}
R~T. Rockafellar and R.~J-B Wets.
\newblock {\em Variational analysis}, volume 317.
\newblock Springer Berlin, Heidelberg, 1998.

\bibitem{shen2024preconditioned}
Xinhua Shen, Zaijiu Shang, and Hongpeng Sun.
\newblock A preconditioned second-order convex splitting algorithm with a difference of varying convex functions and line search.
\newblock {\em arXiv preprint arXiv:2411.07661}, 2024.

\bibitem{Shensun2023}
Xinhua Shen, Hongpeng Sun, and Xuecheng Tai.
\newblock Preconditioned algorithm for difference of convex functions with applications to graph ginzburg–landau model.
\newblock {\em Multiscale Modeling \& Simulation}, 21(4):1667--1689, 2023.

\bibitem{Wen2018}
Bo~Wen, Xiaojun Chen, and Ting~Kei Pong.
\newblock A proximal difference-of-convex algorithm with extrapolation.
\newblock {\em Comput. Optim. Appl.}, 69(2):297--324, Mar 2018.

\bibitem{yang2024proximal}
Lei Yang.
\newblock Proximal gradient method with extrapolation and line search for a class of non-convex and non-smooth problems.
\newblock {\em J. Optim. Theory Appl.}, 200(1):68--103, 2024.

\bibitem{ZhangHongchao}
Hongchao Zhang and William~W. Hager.
\newblock A nonmonotone line search technique and its application to unconstrained optimization.
\newblock {\em SIAM J. Optim.}, 14(4):1043--1056, 2004.

\end{thebibliography}

\end{document}